\documentclass[a4paper,11pt]{article}
\usepackage[english]{babel}
\usepackage[matrix,arrow,tips,curve]{xy}
\usepackage{amsmath,amssymb,amsthm}
\usepackage{enumerate}

\newcounter{subsubsec}

\input{psfig.tex}
\input{psfig.sty}

\newtheorem{teo}{Theorem}[subsubsec]
\newcommand{\pr}[1]{\ensuremath{\mathbb{P}^{#1}}}

\newcommand{\Supp}{\operatorname{Supp}}
\newcommand{\ldo}{.\,.\,.\,}

\newcommand{\Star}{\operatorname{Star}}
\newcommand{\ph}{\ensuremath{\varphi}}
\newcommand{\lis}[1]{\ensuremath{\{x_1,\ldo,x_{#1}}\}}
\newcommand{\som}[1]{\ensuremath{x_1+\cdots+x_{#1}}}
\newcommand{\fx}{\ensuremath{\Sigma_X}}
\newcommand{\fy}{\ensuremath{\Sigma_Y}}
\newcommand{\Z}{\ensuremath{\mathbb{Z}}}
\newcommand{\f}{\ensuremath{\Sigma}}
\newcommand{\Q}{\ensuremath{\mathbb{Q}}}
\newcommand{\RelInt}{\operatorname{Rel\,Int}}
\newcommand{\Int}{\operatorname{Int}}
\newcommand{\NE}{\operatorname{NE}}
\newcommand{\N}{\ensuremath{\mathcal{N}_1}}
\newcommand{\A}{\ensuremath{\mathcal{A}_1}}
\newcommand{\ld}{.\,.\,}

\newtheorem{cor}[teo]{Corollary}
\newtheorem{emdefi}[teo]{Definition}
\newtheorem{prop}[teo]{Proposition}
\newtheorem{lemma}[teo]{Lemma}
\newtheorem{conj}[teo]{Conjecture}
\newenvironment{dimo}[1][Proof.] {\begin{proof}[#1]} {\end{proof}}
\newcommand{\PC}{\operatorname{PC}}
\newenvironment{defi}{\begin{emdefi} \em}{\end{emdefi}}

\begin{document}

\begin{center}
{\em Toric Fano varieties and birational morphisms \\

\medskip

 Cinzia Casagrande }




\end{center}

\bigskip
\bigskip

\noindent Smooth toric Fano varieties are classified up to dimension 4.
In dimension 2 there are five toric Del Pezzo surfaces: $\pr{2}$,
$\pr{1}\times\pr{1}$ and $S_i$, the blow-up of $\pr{2}$ in $i$ points, for 
$i=1,2,3$. There are 18 toric Fano
3-folds~\cite{bat3,wat} and 124 toric Fano 4-folds~\cite{bat2,sato}. 
In higher dimensions, little is known about them: many
properties that hold in low dimensions, are not known to hold in
general. 

Let $X$ be a toric Fano variety of dimension $n$. 
We recall that in $X$ linear and algebraic equivalence 
for divisors coincide, so 
the Picard number $\rho_X$ is the rank of the
Picard group of $X$. 
By means of some basic combinatorial properties of $\fx$, we show
that \emph{for every irreducible invariant divisor $D\subset X$, we have
$0\leq\rho_X-\rho_D\leq 3$}. Moreover,
\emph{if $\rho_X-\rho_D=3$, then $X$ is a toric
$S_3$-bundle over a lower-dimensional toric Fano variety;
if $1\leq \rho_X-\rho_D\leq 2$ (with some additional hypotheses in the
case $\rho_X-\rho_D=1$), we give an explicit birational
description of $X$}. 
 This gives a structure theorem for a large class of toric Fano
varieties (theorem~\ref{structure}).

The best known bound for the Picard number of $X$ is due to
O.~Debarre~\cite{debarre}: $\rho_X\leq 2\sqrt{2n^3}$. 
On the other hand, we know that the maximal Picard number in dimension $n$
is at least $2n$ for $n$ even, $2n-1$ for $n$ odd. In fact,
up to dimension 4, this is exactly the maximal Picard number.
We show that the same is true in  dimension 5:
\emph{if $X$ is a toric Fano 5-fold, $\rho_X\leq 9$} (theorem~\ref{dim5}).

In the second part of the paper, we study equivariant birational
morphisms $f\colon X\rightarrow Y$ whose source $X$ is Fano. 
We give some combinatorial
conditions on the possible positions of new generators in $\fx$, and
we show that
\emph{for every irreducible invariant divisor $D\subset X$, we have
$\rho_Y-\rho_{f(D)}\leq 3$} (proposition~\ref{malena}).

Then we consider two different special cases of this situation: 
the case where $f$ is a
blow-up, and the case where the dimension is 4.

In the case of $f$ a smooth equivariant blow-up, we are particularly
interested in studying under which hypotheses $Y$ can be
non-projective. In fact, up to now
it is not known whether a toric complete, smooth, non-projective variety
$Y$ can become Fano after a smooth equivariant blow-up 
(while this happens in the non-toric context: see~\cite{bontak},
example~1.3). 
We show that \emph{if $Y$ is non-projective,
 the center $A\subset
Y$ of the blow-up must have dimension $\dim A\geq 3$ and can not be a
projective space} (theorem~\ref{sushi}). 
Hence the first possible case of such a
behaviour should be in dimension 6, with $\dim A=3$.

Finally, we consider the case of $f$ an equivariant birational
morphism $f\colon X\rightarrow Y$ with $X$ Fano and
$\dim X=\dim Y=4$.  We prove a factorization result which was 
announced in~\cite{cras}: \emph{every equivariant birational morphism
between toric 4-folds, whose source is Fano, is a composite of smooth
equivariant blow-ups between toric 4-folds} (theorem~\ref{facto}).
This result is obtained through the  complete
classification of the possible subdivisions in $\fx$ of a maximal cone
$\sigma\in\fy$ (proposition \ref{subdiv}).

The approach we use here to study fans of toric varieties is
the language of primitive collections and primitive relations,
which was introduced by V.~Batyrev. This approach is very useful to deal
with toric Mori theory and with fans of toric Fano varieties: it led to 
Batyrev's classification in dimension 4 \cite{bat2}, and was also used by 
H.~Sato in \cite{sato} in the same context, giving some interesting 
factorization properties of equivariant birational maps and morphisms
between toric Fano 3-folds and 4-folds. 

\medskip

\noindent\emph{Acknowledgments}.
This work is part of my Ph.D.\ thesis. I wish to thank my advisor,
Lucia Caporaso, for constantly guiding me, and 
 for all she taught me in these years.
  
I am also very grateful to Laurent Bonavero, who first pointed out to me
Sato's paper on toric Fano varieties~\cite{sato}, 
for many helpful discussions and hints, and for the kind hospitality
at the Institut Fourier of Grenoble.
 
Finally, I thank Olivier Debarre for some valuable conversations.

\bigskip

{\bf \noindent Contents: }

\medskip

\noindent 1. Preliminaries: primitive relations and toric Mori theory
\dotfill\pageref{prima}\\
2. Primitive collections of order two and invariant divisors
\dotfill\pageref{seconda}\\ 
3. Picard number of toric Fano varieties\dotfill\pageref{terza}\\
4. Equivariant birational morphisms whose source is
Fano\dotfill\pageref{quarta}\\
5. Equivariant birational morphisms in dimension 4: analysis of
the possible subdivisions\dotfill\pageref{quinta}

\medskip

\noindent References\dotfill\pageref{ultima}

\bigskip

\stepcounter{subsubsec}
\subsubsection{Preliminaries: primitive relations and toric Mori theory}
\label{prima} In this section, we recall 
the basic definitions and properties of primitive collections and
primitive relations, and some results of toric Mori theory.
For a more detailed account, we refer to \cite{fulton, oda} for
properties of toric 
varieties, to \cite{bat1,bat2,sato} for primitive collections and
primitive relations and to \cite{reid} for toric Mori theory.

An $n$-dimensional toric variety $X$ is described by a 
finite fan $\fx$ in the
vector space $N_{\Q}=N\otimes_{\Z}\Q$, where $N$ is a free abelian group of
rank $n$. Throughout
the paper we will deal with toric varieties that are smooth and
complete;  this is equivalent to ask that the support of
$\fx$ is the whole space $N_{\Q}$ and every cone in $\fx$ is generated
by a  part of a basis of $N$.

We recall that for each $r=0,\ldo,n$ there is a bijection between
the cones of dimension $r$ in $\fx$ and the orbits of codimension $r$
in $X$; we'll denote by $V(\sigma)$ the closure of the orbit
corresponding to $\sigma\in\fx$ and $V(x)=V(\langle x\rangle)$ in case
of 1-dimensional cones.

For each 1-dimensional cone $\rho\in\fx$, let $v_{\rho}\in\rho\cap N$
be its primitive generator, and 
\[  G(\fx)=\{v_{\rho}\,|\,\rho\in\fx\}  \]
the set of all generators of $\fx$.

The following definitions were introduced by V.~Batyrev in  \cite{bat1}:
\begin{defi}
A subset $P=\lis{h}\subseteq G(\fx)$ is a \emph{primitive collection}
for $\fx$ if $\langle x_1,\ldo,x_h\rangle
\notin\fx$, and $\langle x_1,\ld,\check{x}_i,\ld,x_h\rangle\in\fx$
for each $i=1,\ldo,h$.
 We denote by $\PC(\fx)$ the set of all primitive collections 
of $\fx$.  
\end{defi}
\begin{defi}
Let $P=\lis{h}\subseteq G(\fx)$ be a primitive collection and let 
$\sigma_P=\langle y_1,\ldo,y_k \rangle$ be the unique cone in $\fx$ such that
$\som{h}\in \RelInt \sigma_P$.
Then we get a linear relation
\[
\som{h}-(a_1y_1+\cdots+ a_ky_k)=0
\]
with $a_i$ a positive integer for each $i=1,\ldo,k$. This is the
\emph{primitive relation} associated to $P$. The \emph{degree} of $P$
is the integer $\deg P=h-\sum_{i=1}^k a_i$.
\end{defi}

Let $\A(X)$ be the group of algebraic 1-cycles on $X$ modulo numerical 
equivalence, $\N(X)=\A(X)\otimes_{\Z}\Q$, and $\NE(X)\subset \N(X)$ the
cone of Mori, generated by classes of effective curves.
We recall that there is an exact sequence:
\begin{equation}
\label{sumeri}
\tag{$\heartsuit$}
0\longrightarrow \A(X)\longrightarrow
\mathbb{Z}^{G(\fx)}\longrightarrow N \longrightarrow 0,
\end{equation}
where the map $\A(X)\rightarrow \mathbb{Z}^{G(\fx)}$ is given by $\gamma\mapsto \{\gamma\cdot V(x)\} _{x\in G(\fx)}$. Hence the group $\A(X)$ is canonically isomorphic to the lattice of
integral relations among the elements of $G(\fx)$: a relation
\[ \sum_{x\in G(\fx)} a_x x=0, \qquad a_x\in\Z, \]
corresponds to a 1-cycle that has intersection $a_x$ with $V(x)$ for
all $x\in G(\fx)$.
We will often identify classes in $\N(X)$ and
the associated relations.
We remark that since the canonical class on $X$ is given by
$K_X=-\sum_{x\in G(\fx)}  V(x)$, if $\gamma\in\N(X)$ corresponds to 
$\sum_x a_x x=0$, then $-K_X\cdot\gamma=\sum_x a_x$.
In particular, for every primitive collection $P\in\PC(\fx)$, the
associated primitive relation defines a class $r(P)\in\A(X)$, and
$-K_X\cdot r(P)=\deg P$.

For convenience, we will always write primitive relations as
\[\som{h}=a_1y_1+\cdots+ a_ky_k \]
instead of $\som{h}-(a_1y_1+\cdots+ a_ky_k)=0$, i.~e.\ writing
elements with negative coefficient on the right side. This must not 
be confused with the relation $-(\som{h})+ a_1y_1+\cdots+ a_ky_k=0$,
which is the opposite element in $\A(X)$.

The cone of effective curves in a smooth complete
toric variety has been studied by M.~Reid
in~\cite{reid}, where he shows that in this case $\NE(X)$ is closed
and polyhedral, generated by classes of invariant curves
(\cite{reid}, proposition 1.6). Moreover,  
 he gives a precise description of
the geometry of the fan around a cone corresponding
to an extremal curve:
\begin{teo}[Reid \cite{reid}] 
\label{reidextr}
Suppose that $X$ is projective; let $R$ be an extremal ray of $\NE(X)$ and
 $\gamma$ a primitive element in $R\cap\A(X)$.
 Then the
  relation associated to $\gamma$ has the form
\[ \som{h}=a_1y_1+\cdots+a_ky_k,
\]
with $a_i\in\Z_{>0}$; $P=\lis{h}$ is a primitive collection and
$\gamma=r(P)$. Moreover, for each cone $\langle z_1,\ld,z_t\rangle$ 
such that \[
\{ z_1,\ld,z_t\}\cap\{y_1,\ld,y_k,x_1,\ld,x_h\}=\emptyset\text{ and
}\langle y_1,\ld,y_k,z_1,\ld,z_t \rangle \in\fx,\] 
we have
\[\langle
x_1,\ld,\check{x}_i,\ld,x_h,y_1,\ld,y_k,z_1,\ld,z_t \rangle
\in\fx\quad\text{ for all } i=1,\ldo,h.
\]
\end{teo}
\noindent As 
a consequence of this result, we have an important description of
the cone of effective curves for projective toric varieties:
\begin{prop}[Batyrev \cite{bat1}]
Suppose that $X$ is projective: then
the cone of effective curves $\NE(X)$ is generated by primitive relations.
\end{prop}
\begin{dimo}
It is an easy consequence of 
theorem \ref{reidextr} and lemma \ref{criterio} below.
\end{dimo}
\noindent This 
gives  the following simple characterization of toric Fano varieties:
\begin{prop}[Batyrev \cite{bat2}]
\label{parigi}
Suppose that
$X$ is projective. Then $X$ is Fano if and only if all primitive
collections in $\fx$ have strictly positive degree.
\end{prop}
\noindent Finally, we recall the fundamental result:
\begin{prop}
\label{gradouno}
Let $X$ be a toric Fano variety and $\gamma\in\NE(X)$.
If $\gamma$ has anticanonical degree one, then it is extremal.
\end{prop}
\noindent We 
notice that this result seems well-known by experts in toric geometry
(see V.~Batyrev~\cite{bat2}, theorem~2.3.3 and corollary~2.3.4, and
H.~Sato~\cite{sato},  the proof of lemma~5.4), 
but we couldn't find an explicit proof of it written anywhere.
Hence, we refer to \cite{contr}, proposition 4.3. 

\stepcounter{subsubsec}
\subsubsection{Primitive collections of order two and invariant divisors}
\label{seconda} Let $X$ be a toric Fano variety of dimension $n$.
All toric varieties are birationally equivalent and we know, after the 
 the weak factorization theorem \cite{morelli,abr1}, that 
any two smooth, complete toric varieties can be obtained one 
from the other by a sequence of smooth equivariant blow-ups and blow-downs.

Here we give an \emph{explicit} birational description of $X$
when $\fx$  has at least one pair of symmetric generators $x,-x$. 
This will be linked to some properties of irreducible invariant divisors 
on $X$.

Our starting point is the following lemma:
\begin{lemma}
\label{idea}
Suppose that $\{x,-x\}$ is a
primitive collection in $\fx$. Then for any other primitive collection $P$
containing $x$, the associated primitive relation is  
\[r(P):\quad x+y_1+\cdots+y_h=z_1+\cdots+z_h; \] 
moreover, $P'=\{-x,z_1,\ldo,z_h\}$ is also a primitive collection, with relation 
\[r(P'):\quad  (-x)+z_1+\cdots+z_h=y_1+\cdots+y_h. \]
Both $P$ and $P'$ have degree 1, they are extremal, and $2h\leq n$.
Furthermore, for any other primitive collection
$Q=\{x,u_1,\ldo,u_m\}$ different from $P$ and from $\{x,-x\}$,
we have $V(\langle u_1,\ld,u_m\rangle)\cap V(\langle
z_1,\ld,z_h\rangle)=V(\langle u_1,\ld,u_m\rangle)\cap V(\langle
y_1,\ld,y_h\rangle)=\emptyset$.
\end{lemma}
\noindent For the proof of lemma \ref{idea}, 
we need to recall the following condition for 
effectiveness for a numerical class (see \cite{contr}, lemma 1.4):
\begin{lemma}
\label{criterio}
Let $X$ be a smooth, complete toric variety and
$\gamma\in\A(X)$ given by the relation
\[ a_1x_1+\cdots+a_hx_h-(b_1y_1+\cdots+b_ky_k)=0
\]
with $a_i,b_j\in\Z_{>0}$ for each $i,j$. If $\langle
y_1,\ldo,y_k\rangle\in\fx$, then $\gamma\in\NE(X)$.
\end{lemma}
\begin{dimo}[Proof of lemma \ref{idea}]
Suppose that the primitive relation associated to $P$ is 
\[ x+y_1+\cdots+y_h=a_1 z_1+\cdots+a_k z_k. \]
Since $X$ is Fano, $1\leq \deg P=1+h-\sum_i a_i$, thus $\sum_i a_i \leq h$.
On the other hand, consider the class $\gamma$ in $\N(X)$
corresponding to the relation: 
\[ (-x)+a_1 z_1+\cdots+a_k z_k-(y_1+\cdots+y_h)=0. \]
Since $\langle y_1,\ldo,y_h\rangle\in \fx$, lemma~\ref{criterio}
implies that $\gamma\in\NE(X)$, hence $1\leq -K_X\cdot \gamma=1+\sum_i
a_i-h$, namely  
$\sum_i a_i\geq h$: therefore we get $\sum_i a_i=h$ and $\deg
P=\gamma\cdot (-K_X)=1$. Then $r(P)$ and $\gamma$ are extremal: in
particular, $\gamma$ is primitive, which means $a_i=1$ for all $i$ and
$k=h$. 

Suppose now that $Q=\{x,u_1,\ldo,u_m\}$ is a primitive collection 
different from $P$ and from $\{x,-x\}$.
We want to show that $\langle u_1,\ld,u_m,z_1,\ld,z_h\rangle$ and
$\langle u_1,\ld,u_m,y_1,\ld,y_h\rangle$ are not in $\fx$.
Suppose that $\langle u_1,\ld,u_m,z_1,\ld,z_h\rangle\in\fx$.
Since $r(P)$ is extremal, by theorem~\ref{reidextr} we get that
either $\langle x,u_1,\ld,y_m\rangle\in\fx$, or $P=Q$, in both cases
a contradiction. Suppose now that 
the cone $\langle u_1,
\ld,u_m,y_1,\ld,y_h\rangle$ is in $\fx$; then since
$r(P')$ is extremal and $-x\not\in Q$, 
theorem~\ref{reidextr} implies that also $\langle 
u_1,\ld,u_m,z_1,\ld,z_h\rangle$ is in
$\fx$, which gives a contradiction.
\end{dimo}

\noindent\textbf{Fibration in $\pr{1}$: basic construction.} 
We want  to explain  
the basic geometric construction which is an immediate consequence of
lemma~\ref{idea}. Even if the precise statement of our result will be given 
in theorem~\ref{structure}, we think it is useful for the reader 
to explain this construction here.

We fix two symmetric generators $x,-x\in G(\fx)$.
If the class $x+(-x)=0$ is extremal, then $X$ is a toric 
$\pr{1}$-bundle over a toric Fano variety $W$ and $V(x)$, $V(-x)$ are the 
invariant sections of the bundle. In this case $x$ and $-x$ are not contained
in any primitive collection different from $\{x,-x\}$.
Otherwise, if there exists a primitive collection $P\neq\{x,-x\}$ 
containing $x$,  the class $x+(-x)=0$ is not extremal:
 a decomposition of this class in $\NE(X)$ is given by $r(P)+r(P')$, 
where $P'$ is as in lemma~\ref{idea}.

For every such $P$, we define   $E_P$ to be the union of the exceptional 
loci of $r(P)$ and $r(P')$. 
By lemma~\ref{idea}, for two different $P$, $Q$ containing $x$, the loci 
$E_P$ and $E_Q$ are disjointed. 
The locus $E_P$ has pure
codimension $\# P-1$, hence the divisorial ones are 
exactly the ones for $P$ of order two.

\smallskip

\noindent\emph{Fact:} outside the loci of the form $E_P$, $X$ is fibered in 
$\pr{1}$. 

\smallskip

\noindent We think of the loci $E_P$ as the ``obstructions'' to $X$ being 
a $\pr{1}$-bundle.
We can ``eliminate'' these obstructions by a finite number of flips
and blow-downs, getting a smooth projective toric variety which is 
a toric $\pr{1}$-bundle over $W\simeq V(x)$.
More precisely:

\noindent $\bullet\ $ for every primitive collection $P=\{x,y\}$
of order two, we blow-down 
the divisor $V(y)$ with respect to the extremal relation $r(P'):\ 
(-x)+z=y$;

\noindent $\bullet\ $ 
for every primitive collection $P=\{x,y_1,\ld,y_h\}$ with $h\geq 2$,
we flip the relation $r(P'):\ (-x)+z_1+\cdots+z_h=y_1+\cdots+y_h$:
 namely, we blow-up the subvariety
$V(\langle y_1,\ldo,y_h \rangle)\subset E_P$ 
and then blow-down the exceptional
divisor $V(v)$   with respect to the extremal class
$(-x)+z_1+\cdots+z_h=v$; in this way we get a new smooth projective
(non-Fano) toric variety with an extremal class given by
the primitive relation $y_1+\cdots+y_h=(-x)+z_1+\cdots+z_h$.
 The flip exchanges an extremal
class of anticanonical degree 1 with an extremal class of
anticanonical degree -1.

We perform this for all primitive collections containing $x$ 
(different from $\{x,-x\}$). In the end we get a smooth projective
toric variety $\overline{X}$ 
which is a toric $\pr{1}$-bundle over a smooth toric 
variety $W$. The images of the divisors $V(x)$, $V(-x)$ in $\overline{X}$ 
are the invariant sections of the bundle;
since the divisor $V(x)$ in $X$ is disjoint from the exceptional 
loci of the flips and of the blow-downs, we have $V(x)\simeq W$.

\smallskip

\noindent\emph{Geometric description in dimension 4:}
in dimension 4, when $E_P$ is not divisorial, it is the union of two 
$\pr{2}$ which intersect in one point. Each $\pr{2}$  intersects
one (and only one) of the divisors $V(x)$, $V(-x)$ along a $\pr{1}$.
Blowing-up the $\pr{2}$ that touches $V(-x)$, we get an exceptional divisor
isomorphic to $\pr{2}\times\pr{1}$ with normal bundle $\mathcal{O}(-1,-1)$.
The other $\pr{2}$ is blown-up in one point. 
When we blow-down the exceptional divisor
on $\pr{1}$, the image of $E_P$ is  $S_1$, fibered in $\pr{1}$ on the
 intersections with the images of 
the two divisors $V(x)$, $V(-x)$. The divisor $V(-x)$ is ``flopped''.

\vspace{15pt} 

\hspace{10pt}\psfig{figure=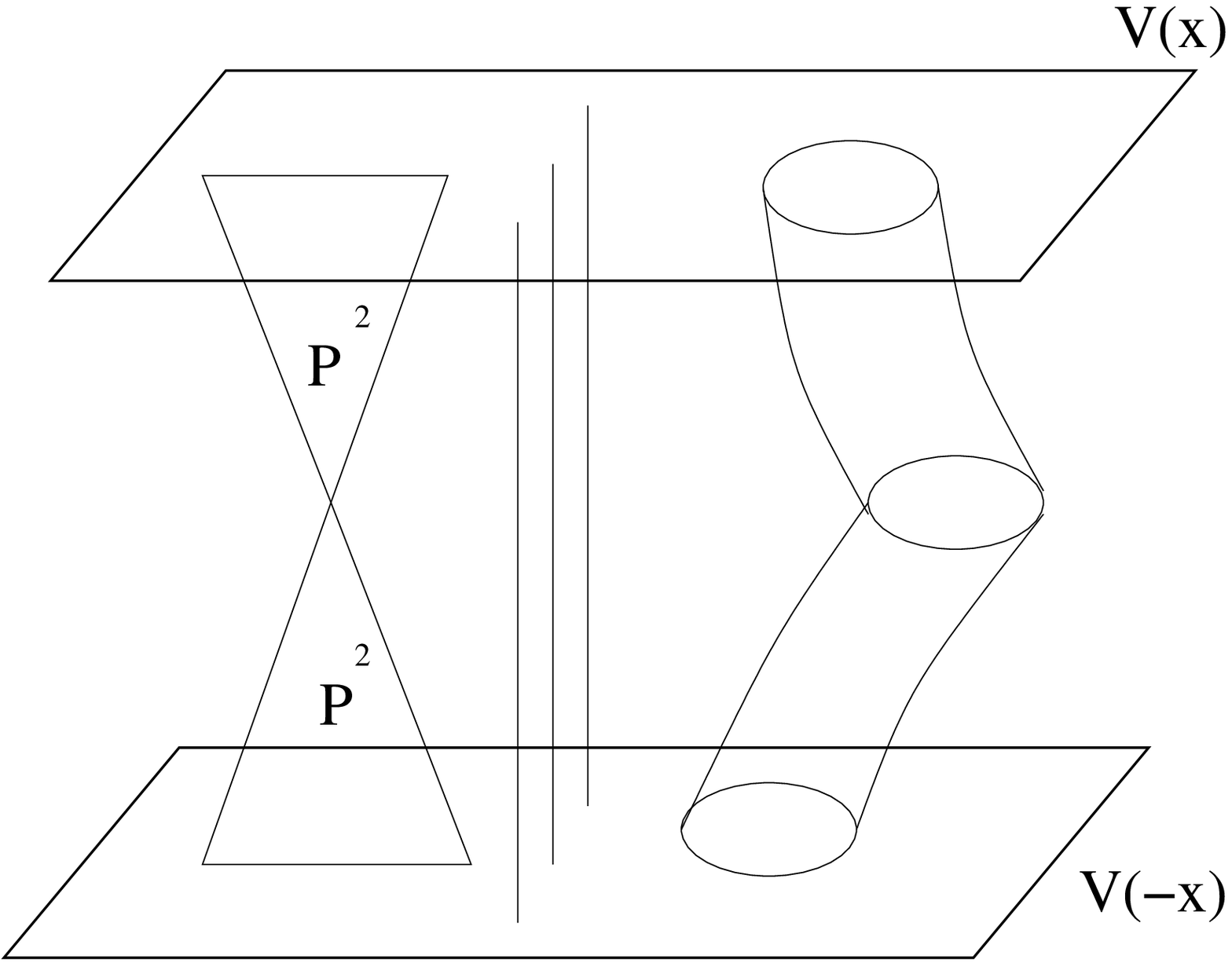,width=4.5cm}
\hspace{40pt}\psfig{figure=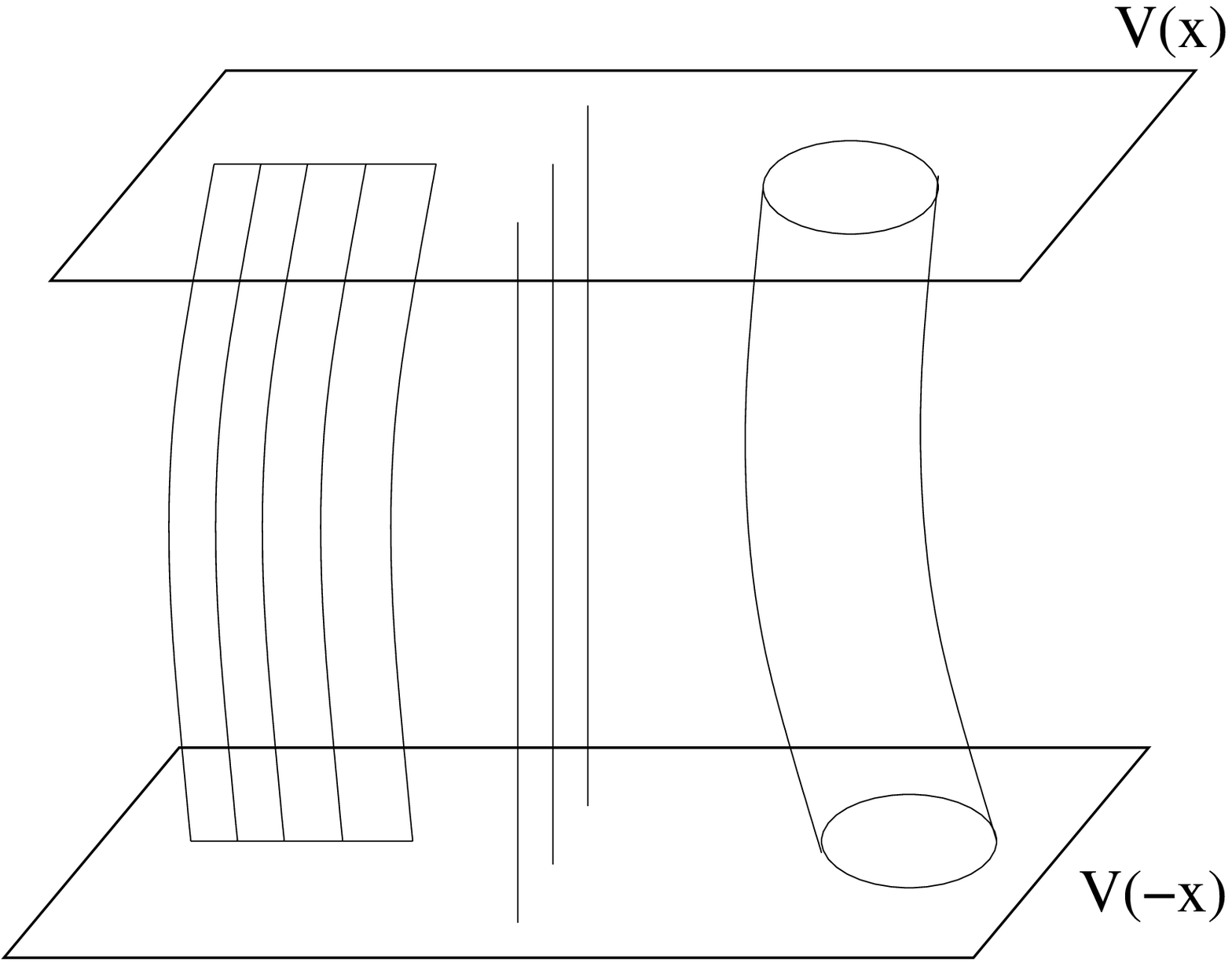,width=4.5cm}
\vspace{15pt}

In general, the loci $V(\langle 
y_1,\ldo,y_h\rangle)$ and $V(\langle z_1,\ldo,z_h\rangle)$ are 
toric $\pr{h}$-bundles over a smooth toric variety $T$ of dimension $n-2h$,
and intersect along a common invariant section.
Fiberwise on $T$, the same description as in the 4-dimensional case holds
(with $\pr{h}$ instead of $\pr{2}$).

\smallskip

\noindent\emph{Remark:} we recall that a toric bundle is locally trivial
on Zariski open subsets; see \cite{ewald2}, theorem 6.7 on page 246.

\medskip

Our next lemma concerns primitive collections of order two in $\fx$.
 Due to the fact that $X$ is Fano, 
 these primitive collections 
must have a particularly simple primitive 
relation.
In fact, if $\{x,y\}$ is a primitive
collection, its degree must be positive by proposition~\ref{parigi},
hence there are only two
possible primitive relations: $x+y=0$ or $x+y=z$. In this last case,
the relation has degree 1, so by proposition~\ref{gradouno} 
the corresponding class in $\NE(X)$ is
extremal. 
Therefore, for each pair of generators $x,y\in G(\fx)$, only  three
cases can occur:  
\begin{enumerate}[(i)]
\item $\langle x,y\rangle\in\fx$
\item $y=-x$ 
\item $z=x+y\in G(\fx)$, and the relation $x+y=z$ is extremal, hence
  by theorem~\ref{reidextr} 
the cones $\langle x,z\rangle$ and $\langle y,z\rangle$ are in  $\fx$. 
\end{enumerate}
This has a nice geometric interpretation. We recall  
 that two irreducible invariant divisors $V(x)$, $V(y)$ in $X$ are
 disjoint if and only if the cone $\langle x, y\rangle$ is not in
 $\fx$. Hence, if $V(x)\cap V(y)=\emptyset$, there are only two
 possibilities:
\begin{enumerate}[(a)]
\item $x+y=0$, hence $X$ is generically fibered in $\pr{1}$ as
  described in the basic construction;
\item there exists a third irreducible invariant divisor $D$ that is a
  toric $\pr{1}$-bundle over a smooth toric variety $Z$; $D\cap V(x)$
  and $D\cap V(y)$ are the sections of the $\pr{1}$-bundle. The normal
  bundle of $D$ restricted to a fiber is $\mathcal{O}_{\pr{1}}(-1)$,
  hence there exists a smooth projective variety $Y$ and a smooth
  equivariant blow-up $X\rightarrow Y$ with exceptional divisor
  $D$. In $Y$ the images of $V(x)$ and $V(y)$ intersect along a
  codimension 2 subvariety which is the center of the blow-up; 
  $Y$ can fail to be Fano.
\end{enumerate}
\begin{lemma}
\label{relations} 
If $\fx$ has two different primitive relations
$x+y=z$ and $x+w=v$, then $w=-x-y$ and $v=-y$.
Therefore there are at most two primitive collections of order 2 and degree
1 containing $x$, and the associated primitive relations are
$x+y=(-w)$ and $x+w=(-y)$.
\end{lemma}
\begin{dimo}
Suppose $y\neq w$. The two relations $x+y=z$ and $x+w=v$ are extremal, 
thus the cones
$\langle x,z\rangle$, $\langle y,z\rangle$, $\langle x,v\rangle$, 
$\langle w,v\rangle$ are in $\fx$. 

If $\langle z,v\rangle\in\fx$, then 
$\langle z,v,y\rangle$ and $\langle z,v,w\rangle$ are in $\fx$, which
is impossible, because $z+w=v+y$. Therefore $\{z,v\}$ is a primitive
collection. 

If the associated primitive relation is $z+v=0$, we have $2x+y+w=0$. 
We claim that this is impossible. Indeed, if $2x+y+w=0$
the set $\{y,w\}$ can not be
a part of a basis of the lattice, because $x$ is not an integral linear
combination of $y$ and $w$. Since  $X$ is smooth, the cone $\langle
y,w\rangle$ is not in $\fx$, so $\{y,w\}$ is a primitive collection.
Its primitive relation can not be $y+w=0$ (it would imply $x=0$)
nor $y+w=t$ (it would be $t=-2x$), so its degree
can not be positive, which is a contradiction because $X$ is Fano.

Therefore the primitive
relation associated to $\{z,v\}$ 
is $z+v=u$.  We want to show that $u=x$, which gives the statement.

If $u=y$, we get $2x+w=0$, which is impossible because the generators
are primitive elements in $N$. The same if $u=w$.

Suppose that $u\not\in\{x,y,w\}$. Since the degree of the relation $z+v=u$
is 1, this class
is extremal:  so $\langle u,v\rangle$
and $\langle u,z\rangle$ are in $\fx$, and since $x+y=z$ and $x+w=v$
are extremal, all the cones
$\langle u,v,x\rangle$, $\langle u,v,w\rangle$, $\langle
u,z,x\rangle$, $\langle u,z,y\rangle$ are in $\fx$.

\vspace{15pt} 

\hspace{100pt}\psfig{figure=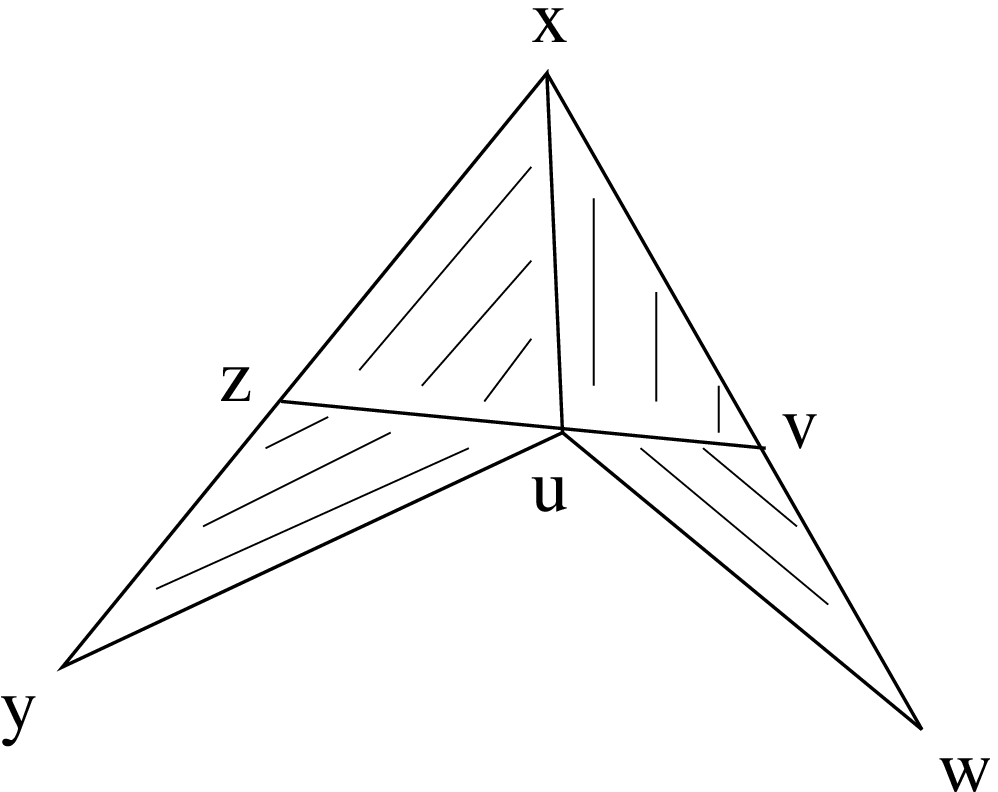,width=4cm}
\vspace{15pt}

\noindent Now consider the pair $\{y,v\}$: since $y+2v=u+w$, the cone $\langle
y,v\rangle$ can not be in $\fx$,  hence $\{y,v\}$ is a primitive
collection. It can not be $y+v=0$, because it would imply $u=x$.
Moreover, there can not be a primitive relation 
$y+v=k$, because this would imply 
$\langle k,v\rangle\in\fx$, which is impossible since $k+v=u+w$.
Therefore we get a contradiction.
\end{dimo}

We are now ready to state the principal result of this section. 
\begin{teo}
\label{structure}
For every irreducible invariant divisor 
$D\subset X$, we have $0\leq \rho_X-\rho_D\leq 3$. 

\medskip

\noindent \textbf{\ref{structure}.1. Case $\rho_X-\rho_D=1$ and
$D=V(x)$ with $x,-x\in G(\fx)$.}  

\smallskip

\noindent There
   exists a sequence 
\begin{equation}
\label{dash}
\tag{$\diamondsuit$}
X=X_1 \stackrel{\ph_1}{\dasharrow} X_2
\stackrel{\ph_2}{\dasharrow}\cdots
\stackrel{\ph_r}{\dasharrow} X_{r+1}\stackrel{\psi}{\rightarrow}\overline{X}
\end{equation}
where for all $i$, $X_i$ is a smooth projective
toric variety and $\ph_i$ is a
toric flip that exchanges an extremal class having anticanonical
degree 1 in $\NE(X_i)$ with an extremal class having anticanonical
degree -1 in $\NE(X_{i+1})$.

In    (\ref{dash}) $\psi$ is an isomorphism and
$\overline{X}=X_{r+1}$ is a toric
 $\pr{1}$-bundle over a smooth toric variety $W\simeq D$. 
The exceptional loci of the flips are disjoint from $D$, and the
 image of $D$ in $\overline{X}$ is a section of the $\pr{1}$-bundle.
 Moreover, $D$ is Fano if and only if $r=0$.

\medskip

\noindent \textbf{\ref{structure}.2. Case $\rho_X-\rho_D= 2$.} 

\smallskip

\noindent There exists a sequence as (\ref{dash}) such that
one of the following occurs:

\smallskip

\noindent \ref{structure}.2a.
the morphism $\psi$ is an isomorphism and
$\overline{X}=X_{r+1}$ is a toric $S_2$-bundle over a
smooth toric Fano variety $Z$. 

There is an irreducible invariant Fano divisor $D'\subset X$ with
$\rho_{D'}=\rho_D$ such that the exceptional loci of
the flips are disjoint from $D'$. The image $\overline{D}'$
of $D'$ in $\overline{X}$ is a
toric $\pr{1}$-bundle over $Z$:
\[
\xymatrix{
{\overline{D}'\ } \ar[dr]_{\pi_{|\overline{D}'}} \ar @{^{(}->}[r]  &
{\overline{X}} 
\ar[d]^{\pi}  \\ 
 & Z }
\]
 and for all $z\in Z$ $\overline{D}'\cap\pi^{-1}(z)$ is
the exceptional curve in $S_2\simeq \pi^{-1}(z)$
given by the proper transform of the line in $\pr{2}$ through
the two blown-up points. 

\smallskip

\noindent \ref{structure}.2b.  the variety $\overline{X}$ is a toric
 $\pr{1}$-bundle over a smooth toric variety $W\simeq D$ and
$\psi$ is a smooth equivariant 
blow-up of a  codimension 2 subvariety of $\overline{X}$. 

The exceptional loci
 of the flips and of the blow-down are disjoint from $D$, and the
 image of $D$ in $\overline{X}$ is a section of the $\pr{1}$-bundle. 
Moreover, $D$ is Fano if and only if $r=0$.

\medskip

\noindent \textbf{\ref{structure}.3.  Case $\rho_X-\rho_D= 3$. }

\smallskip

\noindent $X$ is a toric $S_3$-bundle
over  a smooth  toric Fano variety $Z$.
 
There are six invariant Fano divisors 
$D_1=D,D_2,\ldo,D_6$ in $X$, all
having Picard number $\rho_X-3$, that are toric $\pr{1}$-bundles over
$Z$:
\[
\xymatrix{
{D_i\ } \ar[dr]_{\pi_{|D_i}} \ar @{^{(}->}[r]  &  X \ar[d]^{\pi}  \\
 & Z }
\]
For all $z\in Z$,  $D_i\cap\pi^{-1}(z)$ are the six
 exceptional curves in $S_3\simeq \pi^{-1}(z)$.
\end{teo}
\noindent\emph{Remark:} 
 an irreducible invariant divisor
$D=V(x)\subset X$ has  always $-K_D$ nef; $D$ is Fano if and only if 
$x$ is not contained in any primitive collection of degree 1 and order
greater than 2, namely if
and only if there are no curves $C\subset D$ 
such that $-K_X\cdot C=D \cdot C=1$ (see ~\cite{bat2}, proposition
2.4.4 and corollary 2.4.5). 
We are interested in finding Fano divisors in $X$ because this would allow 
to show some properties of $X$ using induction. We will come back to this 
in the next section.

\medskip

Tables 1 and 2 on page \pageref{tav}
describe the situation in dimensions 3 and 4.
The number in every box is the number of varieties verifying the 
corresponding case, and the number inside the parentheses in the first row 
and in the left column is the total number of varieties verifying that case.

For the 3-dimensional case, we refer to \cite{oda}, pages 90-92. 
We call $F$ the
toric Fano 3-fold with Picard number 5
which is not a product; it is obtained from 
$\mathbb{P}_{\pr{1}}
(\mathcal{O}\oplus\mathcal{O}\oplus\mathcal{O}(1))$ blowing up the
three invariant sections, and  is a toric $S_3$-bundle on $\pr{1}$.
We see in table 1 that the situation is very simple: the 18 toric Fano
3-folds are all either a toric bundle with fiber
$\pr{1}$, $\pr{2}$, $\pr{3}$, $S_2$ or $S_3$, or they are obtained from a
$\pr{1}$-bundle with a single blow-up (case \ref{structure}.2b, r=0).

In table 2, the letters refer to Batyrev's table of toric Fano 4-folds
in \cite{bat2}. The variety missing in Batyrev's table is Sato's
4-fold, see \cite{sato}. 
We remark first that with only two
exceptions, in all cases described by theorem \ref{structure} there is
at most one flip. The two exceptions are the Del Pezzo variety $V^4$
and the pseudo-Del Pezzo variety $\widetilde{V}^4$; these two
varieties are special from many points of view: see
\cite{VK,ewald,sato,debarre}. 

We remark also that theorem \ref{structure} describes 107 of the 124
toric Fano 4-folds. Among the remaining 4-folds, 9 are projective bundles;
the classes for which we really lack a description are the last two.

\renewcommand{\arraystretch}{1.2}

\begin{table}
\label{tav}
\scriptsize
\[
\begin{array}{|l|c|c|c|c|c|}
\cline{2-6} 
\multicolumn{1}{c|}{}& \rho=1 &\rho=2&\rho=3&\rho=4&\rho=5\\
\multicolumn{1}{c|}{}& (1) & (4)&(7)&(4)&(2)\\
\hline \ref{structure}.1,r=0&&3&5&&\\
\hline
\ref{structure}.2a,r=0&&&&4&\\
\hline\ref{structure}.2b,r=0&&&2&&\\
\hline\ref{structure}.3&&&&&S_3\times 
\pr{1},F\\
\hline\pr{s}\text{-bundles},s>1&\pr{3}&1&&&\\
\hline \end{array}
\]
\caption{toric Fano 3-folds.}

\vspace{35pt}

\[
\begin{array}{|l|c|c|c|c|c|c|c|c|}
\cline{2-9} 
\multicolumn{1}{c|}{}&\rho=1&\rho=2&\rho=3&\rho=4&\rho=5&\rho=6&\rho=7
&\rho=8\\
\multicolumn{1}{c|}{}&(1)&(9)&(28)&(47)&(27)&(10)&(1)&(1)\\
\hline 
\ref{structure}.1,r=0    &&4  &16 &16   &&&&\\
\multicolumn{1}{|r|}{(36)}&&(B)&(D)&(I,L)&&&&\\ \hline
\ref{structure}.1,r=1    &&&2        &4  &&&&\\
\multicolumn{1}{|r|}{(6)}&&&(G_1,G_3)&(M)&&&&\\ \hline
\ref{structure}.1,r=3    &&&&&1&&&\\
\multicolumn{1}{|r|}{(1)}&&&&&(\widetilde{V}^4)&&&\\ \hline
\ref{structure}.1,r=6    &&&&&&1&&\\
\multicolumn{1}{|r|}{(1)}&&&&&&(V^4)&&\\ \hline
\ref{structure}.2a,r=0    &&&&10 &17 &S_2\times S_2&&\\
\multicolumn{1}{|r|}{(28)}&&&&(H)&(Q)&&&\\ \hline
\ref{structure}.2a,r=1   &&&&&3&&&\\
\multicolumn{1}{|r|}{(3)}&&&&&(R)&&&\\ \hline
\ref{structure}.2b,r=0    &&&3  &12 &1      &&&\\ 
\multicolumn{1}{|r|}{(16)}&&&(E)&(I)&(3.4.1)&&&\\ \hline
\ref{structure}.2b,r=1   &&&&2  &&&&\\ 
\multicolumn{1}{|r|}{(2)}&&&&(J)&&&&\\ \hline
\ref{structure}.3         &&&&&4  &8  &S_3\times S_2& S_3\times S_3 \\
\multicolumn{1}{|r|}{(14)}&&&&&(K)&(U)&&\\ \hline
\pr{s}\text{-bundles},        &\pr{4}&5    &3  &&&&&\\
s>1\qquad\quad (9)&      &(B,C)&(D)&&&&&\\ \hline 
\text{no symm.\ gen.}    &&&4  &1    &1&&&\\
\multicolumn{1}{|r|}{(6)}&&&(G)&(M_5)&\text{(Sato's)}&&&\\ \hline
\text{no prim.\ coll.}&&&&2&&&&\\
\text{of order 2 }\ \ (2) &&&&(Z)&&&&\\ \hline
\end{array}
\]
\caption{toric Fano 4-folds.}
\end{table}
\begin{dimo}[Proof of theorem \ref{structure}]
Let $D=V(x)$. We recall that from the exact sequence (\ref{sumeri}) on page 
\pageref{sumeri} we
 have $\rho_X= \# G(\fx)-n$ for the Picard number of $X$, and since 
1-dimensional cones in the fan $\Sigma_D$ are given exactly by
2-dimensional cones of $\fx$ containing $x$, we get:
\begin{eqnarray*}
 \rho_D&=&\#\{y\in G(\fx)\,|\,y\neq x \text{ and }\langle y,x
\rangle\in\fx\}-(n-1)\\
&=& \#G(\fx)-\#\{y\in G(\fx)\,|\,\{y,x\}\text{ is a primitive collection}
\}-n \\
&=&\rho_X-\#\{y\in G(\fx)\,|\,\{y,x\}\text{ is a primitive collection}
\}.
\end{eqnarray*}
Hence the difference $\rho_X-\rho_D$ is equal to the number of
primitive collections 
of order two containing $x$. 
Let's show that $x$ is contained in at most three 
primitive collections of order two. By lemma~\ref{relations}, there
are at most two primitive collections of order two and degree 1
containing $x$, and the only primitive collection of order two and
degree 2 containing $x$ can be $\{x,-x\}$. So we have $\rho_X-\rho_D\leq 3$.

\smallskip

For the second part of the theorem,
the idea of the proof is the following. Suppose that we have a pair
$x,-x\in G(\fx)$. Then $X$ is generically fibered in $\pr{1}$ as
described in the basic construction. Moreover, we know by
lemma~\ref{relations} that there are at most two divisorial $E_P$. 
 We are going to show that if there are exactly two, then there are no
 other obstructions; in this case there are three
 pairs of symmetric generators, which give a toric bundle in surfaces 
$S_3$.
If there is only one divisorial obstruction, we have two pairs of
symmetric generators, and the non-divisorial obstructions to the two
pairs are ``compatible'': hence with a finite number of flips we get 
a toric bundle in surfaces $S_2$.

\smallskip

\noindent 1)  Since $\rho_X-\rho_D=1$, $x$ not contained in primitive 
collections of order two different from $\{x,-x\}$. Hence 
the statement follows from the basic 
construction: none of the $E_P$ is divisorial.

\smallskip

\noindent 2) $D=V(x)$ with $x$ contained in two primitive collections
  of order two.
We consider first the case where both
collections have degree 1: then by lemma~\ref{relations}, the
primitive relations are  $x+y=(-w)$, $x+w=(-y)$. Thus, by lemma~\ref{idea},
we know the following primitive relations of $\fx$:
\[  
(-y)+(-w)=x, \qquad y+(-y)=0,\qquad w+(-w)=0.
\]
Moreover, we know that $\langle y,w\rangle\in\fx$, because if
$\{y,w\}$ were a primitive collection, the primitive relation should be
$y+w=(-x)$, and this would give a third primitive collection $\{x,-x\}$
of order two containing $x$. Hence  in  the plane $H\subset N_{\Q}$
spanned  by $x$ and $y$ we get a fan of the Del Pezzo surface $S_2$ as
in the following figure: 

\vspace{15pt} 

\hspace{100pt}\psfig{figure=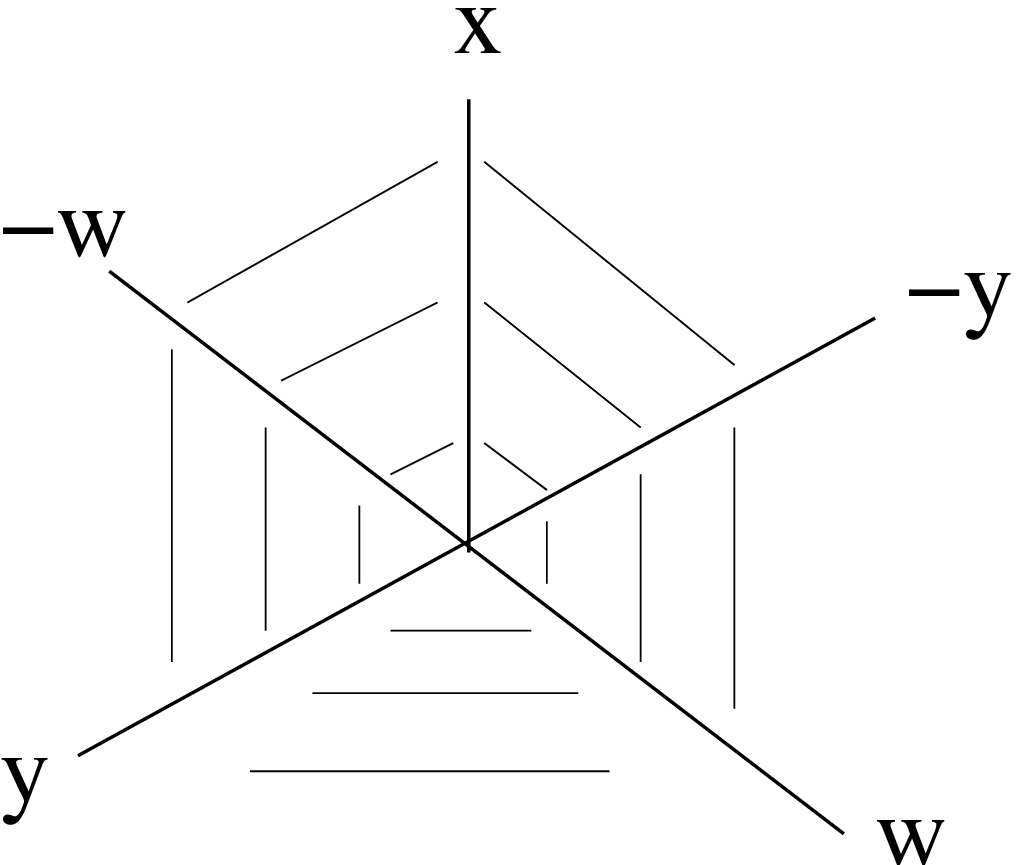,width=3cm}
\vspace{15pt}

\noindent We remark that the generators $x,y,-y,w,-w$ can not be
contained in 
other primitive collections of order two beyond the ones already given.

\smallskip

\noindent\emph{Claim 1:} if $P$ is a primitive collection such that
$\#P\geq 3$ and $y\in P$, then it is also $w\in P$.

\noindent\emph{Proof of claim 1.} Let $P=\{y,z_1,\ldo,z_h\}$. By 
lemma~\ref{relations}, we know that the cone 
$\langle -y,z_1,\ldo,z_h\rangle$ is in $\fx$.
If $w\not\in P$,  the extremality of the relation $(-y)+(-w)=x$
 implies by theorem~\ref{reidextr} that $\langle -w,z_1,\ldo,z_h\rangle\in\fx$;
again, using the relation $x+y=(-w)$, we get $\langle P\rangle\in\fx$,
a contradiction.

\smallskip

Suppose that $P=\{y,w,z_1,\ldo,z_h\}$
with $h\geq 1$. Clearly $\{z_1,\ldo,z_h\}\cap\{x,-y,-w\}=\emptyset$,
otherwise $P$ would contain another primitive collection.
We know by lemma~\ref{idea} that $\deg P=1$ 
and that $\fx$ contains the following extremal relations:
\[
\begin{array}{rcl} 
y+w+z_1+\cdots+z_h&=&v_1+\cdots+v_{h+1}, \\
(-y)+v_1+\cdots+v_{h+1}& =&w+z_1+\cdots+z_h, \\
(-w)+v_1+\cdots+v_{h+1}& =&y+z_1+\cdots+z_h.
\end{array}
\]
Moreover, it is easy to see that also 
\begin{equation}
\label{pelosi}
\tag{$\clubsuit$}
x+v_1+\cdots+v_{h+1} =z_1+\cdots+z_h 
\end{equation}
is a primitive relation; it has degree 2 and it is
not extremal.

\smallskip

\noindent\emph{Claim 2:}  if $P$ is a primitive collection
 such that $\#P\geq 3$ and $x\in P$, then $P$ is 
 obtained as (\ref{pelosi}), so in
particular it has degree 2.

\noindent\emph{Proof of claim 2.}  
Let  $P=\{x,p_1,\ldo,p_r\}$. Then it is
easy to see that also $\{-y,p_1,\ldo,p_r\}$ is primitive, so by
lemma~\ref{idea} we get extremal relations:
\[ (-y)+p_1+\cdots+p_r=q_1+\cdots+q_r, \quad
y+q_1+\cdots+q_r=p_1+\cdots+p_r. \]
By claim 1, we can assume $w=q_r$; then $r(P)$ is 
$x+p_1+\cdots+p_r=q_1+\cdots+q_{r-1}$, which is like 
 (\ref{pelosi}).

\smallskip

Clearly, for symmetry, claim 1 implies that if $P$ is a primitive
collection such that $\#P\geq 3$ and $w\in P$, then it is also $y\in
P$. Now consider the two pairs $y,-y$ and $w,-w$. Applying the basic
construction, by claim 1 we see that we can eliminate 
the obstructions for $y+(-y)=0$ and for $w+(-w)=0$
simultaneously. Namely, we flip all the
relations  of type
$y+w+z_1+\cdots+z_h=v_1+\cdots+v_{h+1}$, and get a smooth toric variety
$\overline{X}$. Consider the blow-down $\overline{X}\rightarrow Y$
with respect to the relation $(-y)+(-w)=x$. In $Y$ the classes
$y+(-y)=0$ and $w+(-w)=0$ are both extremal, thus $Y$ is a toric 
$(\pr{1}\times\pr{1})$-bundle over a smooth toric variety $Z$. The
center of the blow-up $\overline{X}\rightarrow Y$ is an invariant
section, hence $\overline{X}$ is a toric $S_2$-bundle over $Z$. 
Moreover, $D$ is the exceptional divisor of this blow-up, and the
center of the blow-up is isomorphic to $Z$, so
$D$ is a $\pr{1}$-bundle over $Z$. Claim 2 implies that $D$ is Fano,
because  in $\fx$ there are no
primitive collections containing $x$ and having degree 1 and order
greater than 2. 
This implies that $Z$ is Fano
(see~\cite{wisnszurek}; in the toric case it is easy to see this directly).
Hence we have shown \ref{structure}.2a.

\smallskip

We  consider now the case where the two primitive collections
containing $x$ are $\{x,-x\}$ and $\{x,y\}$, with relations $x+(-x)=0$
and $x+y=v$. 

If $-y\in G(\fx)$, then we have primitive relations
$v+(-y)=x$ and $v+(-x)=y$. Hence $v$ is contained in two primitive
collection of order 2 and degree 1, and moreover $-v\not\in G(\fx)$,
because otherwise $\{x,-v\}$ would be a third primitive collection of
order two containing $x$. Then we can consider $v$ instead of
$x$ and get again \ref{structure}.2a. 

In the same way, if $-y\not\in G(\fx)$
but $-v\in G(\fx)$, $y$ is contained in two primitive collections of
order 2 and degree 1, so we get again \ref{structure}.2a.

Therefore we can
assume $-y,-v\not\in G(\fx)$: 
\label{tre}

\vspace{15pt} 

\hspace{100pt}\psfig{figure=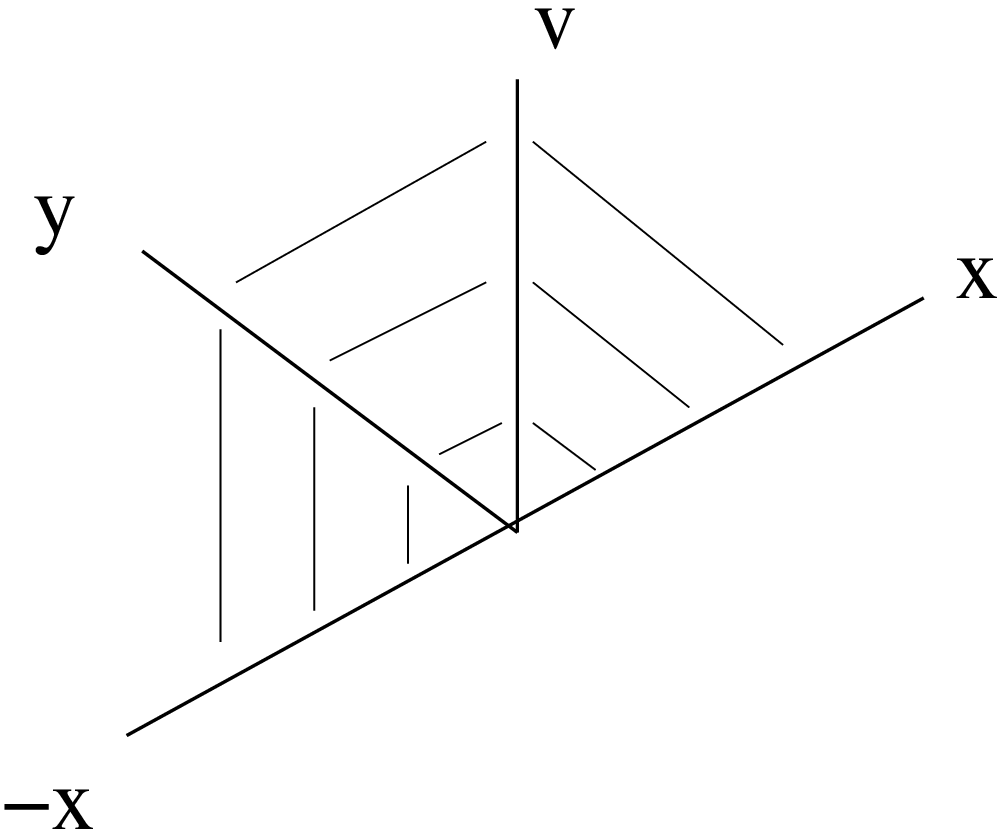,width=3cm}
\vspace{15pt}

\noindent In this case we can perform the basic contruction and get
\ref{structure}.2b; 
we remark that there is one single blow-down corresponding to the relation
$(-x)+v=y$. Moreover, we remark that $D$ is Fano if and only if
there are no other obstructions, namely if and only if $r=0$.

\smallskip

\noindent 3) 
$D=V(x)$ with $x$ contained in three primitive
  collections of order two. 
 By lemma~\ref{relations}, the relations are
\[ x+(-x)=0,\quad x+y=(-w),\quad x+w=(-y).\]
Moreover, applying lemma~\ref{idea} to the pairs $\{x,-x\}$,
$\{y,-y\}$, $\{w,-w\}$
we get the following relations:
\begin{align*} 
 (-x)+(-y)&=w & y+w&=(-x) & y+(-y)&=0 \\ 
 (-x)+(-w)&=y & (-y)+(-w)&=x &w+(-w)&=0.
\end{align*}
The six generators $x$, $y$, $w$, $-x$, $-y$, $-w$ lie in the same plane
 $H\subset N_{\Q}$ as in the following figure:

\vspace{15pt} 

\hspace{100pt}\psfig{figure=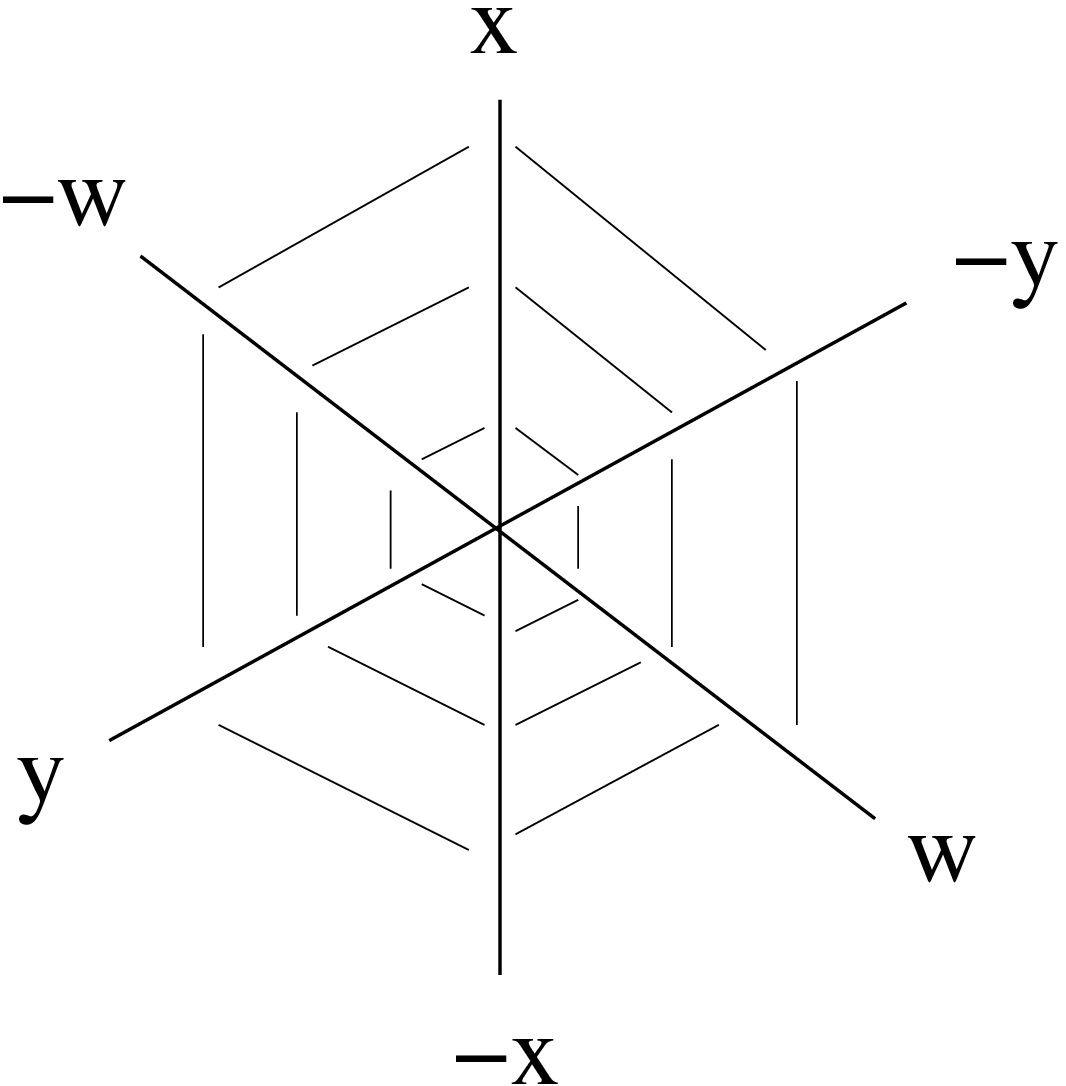,width=3cm}
\vspace{15pt}

\noindent Since all relations having degree
 1 are extremal, the cones $\langle x,-y\rangle$, $\langle
 x,-w\rangle$, $\langle w,-y\rangle$, $\langle -w,y\rangle$, $\langle
 -x,y\rangle$ and $\langle -x,w\rangle$ are all in $\fx$.
These cones give a fan of the Del Pezzo surface $S_3$.

We claim that for every primitive collection $P$ in $\fx$
 different from the nine already given, 
we have $P\cap \{x,y,w,-x,-y,-w\}=\emptyset$.
 Indeed, suppose that $P\cap \{x,y,w,-x,-y,-w\}\neq\emptyset$. 
Since $P$ can not contain another primitive collection, this intersection will have order 1 or 2. If the order is 2, we can
suppose $P=\{x,-y,z_1,\ldo,z_h\}$ with $\{z_1,\ldo,z_h\}\cap
\{x,y,w,-x,-y,-w\}=\emptyset$; but $\langle x,z_1,\ldo,z_h \rangle\in\fx$ 
implies $\langle x,-y,z_1,\ldo,z_h \rangle\in\fx$, a contradiction.
If $P\cap \{x,y,w,-x,-y,-w\}$ has order 1, we can
suppose $P=\{x,z_1,\ldo,z_h\}$ with $\{z_1,\ldo,z_h\}\cap
\{x,y,w,-x,-y,-w\}=\emptyset$. Then by lemma \ref{idea}
we have two primitive extremal relations:
\[ x+z_1+\cdots+z_h=v_1+\cdots+v_h,\quad
 (-x)+v_1+\cdots+v_h=z_1+\cdots+z_h. \]
For what preceeds, we know that $\{v_1,\ldo,v_h\}\cap
\{x,y,w,-x,-y,-w\}=\emptyset$. Then, using theorem \ref{reidextr} and
the primitive extremal relations $(-y)+(-w)=x$, $x+y=(-w)$,
$(-x)+(-w)=y$, we get that $\langle -x,v_1,\ldo,v_h\rangle\in\fx$, a
contradiction. 

Applying the basic construction to the two pairs $x,-x$ and $y,-y$, we
see that \emph{the only obstructions are the divisorial ones}. This
means that blowing-down on $X$ the two relations $x+y=(-w)$ and
$(-x)+(-y)=w$, we get a toric $(\pr{1}\times\pr{1})$-bundle, hence $X$
itself is a toric $S_3$-bundle over a toric Fano variety $Z$.

We remark that the role of the six generators $\{x,y,w,-x,-y,-w\}$
is perfectly symmetric. The six invariant divisors
$D_1=D,D_2,\ldo,D_6$ corresponding to them are all Fano, hence also
$Z$ is Fano.
\end{dimo}

\stepcounter{subsubsec}
\subsubsection{Picard number of toric Fano varieties}
\label{terza} A first bound for the Picard number of $n$-dimensional toric
Fano varieties was due
to V.~E.~Voskresenski\u{\i} and A.~A.~Klyachko~\cite{VK}:
\[ \rho\leq n^2-n+1. \]
Recently, O.\ Debarre~\cite{debarre} has found a better bound:
\[ \rho\leq 2\sqrt{2n^3}. \]
Moreover,  there are some conjectures due to 
G.~Ewald~\cite{ewald} for a 
linear bound in $n$.  

Looking at tables 1 and 2 on page~\pageref{tav}, we see 
that in
dimension 3 and 4 the toric Fano varieties having maximal Picard
number are actually $S_3$-bundles over a toric Fano variety of
lower dimension (and they are very few).
Assuming that this holds in all dimensions, applying
induction would give the following: 
\begin{conj}
\label{bound}
$\rho\leq 2n$  for $n$ even, $\rho\leq 2n-1$ for $n$
  odd. 
\end{conj} 
\noindent This bound would be sharp: just consider 
$(S_3)^{\frac{n}{2}}$ for
$n$ even, $(S_3)^{\frac{n-1}{2}}\times\pr{1}$ or
$(S_3)^{\frac{n-3}{2}}\times F$ for $n$ odd, where $F$ is the toric Fano
3-fold with $\rho_F=5$ which is not a product (see table~1 on
page~\pageref{tav}).  

In this section we show that conjecture~\ref{bound} is true
for $n=5$: 
\begin{teo}
\label{dim5}
Let $X$ be a toric Fano variety of dimension 5. Then $\rho_X\leq
9$. Moreover,  
$\rho_X= 9$ if and only if $X\simeq S_3\times S_3\times\pr{1}$ or 
$X\simeq S_3\times F$.
\end{teo}

Since $\rho_X= \#G(\fx)-n$, giving a bound for the Picard number is
equivalent to giving a bound for the number of generators of $\fx$. 
Let $f_i$ be the number of $(i+1)$-dimensional cones in $\fx$: thus
$f_0$ is the number of generators of the fan and $f_1$ the number of
2-dimensional cones. One has $f_1\leq \binom{f_0}{2}$, and 
$\binom{f_0}{2}-f_1$ is exactly the number of primitive collections of 
order two.

As they count the faces of a simplicial convex polytope, 
the numbers $f_i$ have to satisfy many relations; we refer to 
\cite{mcmullenshephard} for a survey on this. In particular, by the
Dehn-Sommerville equations, $f_{[\frac{n}{2}]},\ldo,f_{n-1}$ are
completely determined by $f_0,\ldo,f_{[\frac{n}{2}]-1}$. In dimension
5 we have (see \cite{mcmullenshephard}, page 112): 
\[
\begin{array}{ccl}
f_2&=&4f_1-10f_0+20, \\
f_3&=&5f_1-15f_0+30, \\
f_4&=&2f_1-6f_0+12.
\end{array}
\]
For toric Fano varieties, we have a further relation, due to V.~Batyrev:
\begin{teo}[\cite{bat2}, 2.3.7]
\label{stino}
Let $X$ be an $n$-dimensional toric Fano variety. Then
\[ 12f_{n-3}\geq (3n-4)f_{n-2}. \]
\end{teo}
\noindent In dimension 5 this gives $12f_2\geq 11 f_3$, that together
with the Dehn-Sommer-ville equations gives: 
\begin{equation}
\label{chao}
\tag{$\spadesuit$}
7f_1\leq 45(f_0-2).
\end{equation}
Using the results of the preceeding sections, 
we can get a new relation under the hypothesis that 
$X$ does not have many primitive collections of order two: 
\begin{lemma}
\label{manu}
Let $X$ be an $n$-dimensional 
toric Fano variety. Suppose that for any generator $x\in
G(\fx)$,  
there are at most two primitive collections of order two containing
$x$, and if there are exactly two, they are $\{x,-x\}$ and $\{x,y\}$ with
$-y,-x-y\not\in G(\fx)$. 
Then $\fx$ has at most $\frac{3}{4}
\,f_0$ primitive collections of order two, namely: 
\[ \binom{f_0}{2}-f_1\leq \frac{3}{4} f_0. \]
\end{lemma}
\begin{dimo}
The maximal number of primitive collections of order two in $X$ is
given by the maximal number of configurations like in
theorem~\ref{structure}.2b, 
where we have three primitive
collections of order two involving four generators
(see figure on page~\pageref{tre}). Hence we can have
$[\frac{f_0}{4}]$ of such configurations, giving  
$3[\frac{f_0}{4}]$ primitive collections of order two. 
It is easy to see that considering also the remaining
$f_0-4[\frac{f_0}{4}]$ generators, the maximal number of primitive
collections of order two is  $\frac{3}{4}f_0$. 
\end{dimo}
\begin{prop}
\label{cland}
Let $X$ be as in lemma~\ref{manu}, and suppose that $\dim X=5$. Then
$\rho_X\leq 8$. 
\end{prop}
\begin{dimo}
By lemma~\ref{manu}, we have 
\[ f_1\geq \frac{1}{2} f_0(f_0-1)-\frac{3}{4} f_0. \]
On the other hand, we have (\ref{chao}): $7f_1\leq 45(f_0-2)$.  These
relations give $14f_0^2-215f_0+360\leq 0$, hence $f_0\leq 13$ and
$\rho_X\leq 8$. 
\end{dimo}
We remark that in this way
 we can not get a similar result in higher dimensions, 
because with theorem~\ref{stino} and the Dehn-Sommerville equations 
we get a relation among  
$f_0,\ldo,f_{[\frac{n}{2}]-1}$. 

\begin{dimo}[Proof of theorem \ref{dim5}]
Suppose that $X$ has three primitive collections of order two
having a common element. Then, by theorem~\ref{structure}, $X$ is a toric
$S_3$-bundle over a toric Fano 3-fold: in particular, $\rho_X\leq 9$.  

Suppose now that $X$ has two primitive collections of order two 
having a common element and that the configuration 
is as described in theorem~\ref{structure}.2a.
Then $X$ has
the same Picard number as a toric $S_2$-bundle over a toric Fano
3-fold; hence $\rho_X\leq 8$. 

Finally, if $X$ does not fit in the preceeding cases, then
proposition~\ref{cland} applies, so $\rho_X\leq 8$. 

Therefore if $\rho_X= 9$, 
 $X$ is a toric
$S_3$-bundle over a toric Fano 3-fold $Z$ with $\rho_Z=4$. 
Looking at table 1 on page \pageref{tav},
 we see that either $Z\simeq S_3\times\pr{1}$ or $Z\simeq F$. 
 In both cases, the only
possibility for $X$ being Fano is $X\simeq S_3\times Z$. 
\end{dimo}

\stepcounter{subsubsec}
\subsubsection{Equivariant birational morphisms whose source is Fano}
\label{quarta} Throughout this section, $X$ will be an $n$-dimensional
toric Fano 
variety, $Y$ an $n$-dimensional, smooth, complete toric variety, and
$f\colon X\rightarrow Y$ an equivariant birational morphism. 
The fan of $X$ is a subdivision of the fan of $Y$, so in particular
$G(\fy)\subseteq G(\fx)$. The new generators of $X$, namely the ones
in $G(\fx)\smallsetminus G(\fy)$, correspond to the irreducible
components of the exceptional divisor of $f$. We want to study which
conditions are imposed on the possible positions of new generators by
the fact that $X$ is Fano.

We fix an irreducible invariant divisor $E=V(x)\subset X$ and consider $A=f(E)=V(\eta)\subset Y$, where $\eta\in\fy$.
Consider $\Star(\eta)=\{\sigma\in\fy|\sigma\supseteq\eta\}$ and
\[
\Lambda_{\eta} =\Supp\Star(\eta)=\bigcup_{\sigma\supseteq\eta} \sigma
\subset N_{\Q}. \]
$\Lambda_{\eta}$ is a closed, connected subset of
$N_{\Q}$, and it has non-empty interior: we think of its boundary as a
partition of the vector space  $N_{\Q}$.
We define the sets of generators:
\[
\mathcal{G}_{\eta}=\{z\in G(\fx)|z\in\Int\Lambda_{\eta}\},\quad
\mathcal{H}_{\eta}=\{z\in G(\fx)|z\not\in\Lambda_{\eta}\}.
\]
$\mathcal{G}_{\eta}$ are the generators of $\fx$ lying inside $\Lambda_{\eta}$,
$\mathcal{H}_{\eta}$ the generators lying outside: we have
$G(\fx)=\mathcal{G}_{\eta}\cup  
\mathcal{H}_{\eta}\cup\{z\in G(\fx)|z\in\partial\Lambda_{\eta}\}$.
\begin{prop}
\label{malena}
In the above setting, we have:
\begin{enumerate}
\item either $\mathcal{H}_{\eta}=\emptyset$, 
or $|\mathcal{G}_{\eta}\cup\mathcal{H}_{\eta}|\leq 4$;
\item for all $z_1\in \mathcal{G}_{\eta}$ and $z_2\in
  \mathcal{H}_{\eta}$, either $z_1+z_2=0$, or $z_1+z_2\in
  G(\fx)\cap\partial \Lambda_{\eta}$; 
\item if $|\mathcal{G}_{\eta}\cup\mathcal{H}_{\eta}|= 4$, then 
$X$ is a toric $S_3$-bundle over a toric Fano variety.
\end{enumerate}
\end{prop}

\medskip

\noindent \emph{Remarks:} 

\medskip

\noindent 1)  if $\mathcal{H}_{\eta}\neq\emptyset$, the majority of
the new generators (the ones in $G(\fx)\smallsetminus G(\fy)$) must
lie on the boundary of $\Lambda_{\eta}$. 

\medskip

\noindent 2) All 
generators in $\mathcal{G}_{\eta}$ are new generators: 
indeed, they forcely lie inside  some
cone of $\fy$. So they all correspond to components of the exceptional
divisor. Instead, the generators in $\mathcal{H}_{\eta}$ can either be
new generators, or they can be generators of $\fy$ which do not belong
to $\Star(\eta)$; such generators are exactly $\rho_Y-\rho_A$, so we
get $\rho_Y-\rho_A\leq |\mathcal{H}_{\eta}|\leq 3$. 

\medskip

\noindent 3) The sets  $\mathcal{G}_{\eta}$  and $\mathcal{H}_{\eta}$
can be seen in a more geometric way: let $z\in G(\fx)$ and
$\tau\in\fy$ such that $z\in\RelInt\tau$. Then $f(V(z))=V(\tau)$ in
$Y$, and we have: 
\[
\begin{array}{ccccc}
z\in\Int \Lambda_{\eta}\ &\Leftrightarrow\ &\tau\supseteq\eta\
&\Leftrightarrow\ & f(V(z))\subseteq A;\\
z\not\in\Lambda_{\eta}\ &\Leftrightarrow\ &\tau+\eta\not\in\fy\ &\Leftrightarrow\ & f(V(z))\cap A=\emptyset.
\end{array}
\]
Therefore we have:
\[
\mathcal{G}_{\eta}=\{z\in G(\fx)|f(V(z))\subseteq A\},\quad
\mathcal{H}_{\eta}=\{z\in G(\fx)|f(V(z))\cap A=\emptyset\}.
\]

\begin{cor}
\label{puh}
Let $X$ be Fano, $Y$ a smooth complete toric variety and 
$f\colon X \rightarrow Y$  an equivariant birational morphism.
Then for any irreducible invariant divisor $E\subset X$, we have 
$\rho_Y-\rho_{f(E)}\leq 3$.

If $f(E)$ is a point  $p\in Y$, then $\rho_Y\leq 3$, and 
 $E$ is the only 
component of the exceptional divisor contracted to $p$, unless
$Y\simeq\pr{n}$. In the case  $Y\simeq\pr{n}$,
there can be at most another component
  contracted to $p$. 
\end{cor}
In particular, if one component of the exceptional
divisor is contracted to a curve, then $\rho_Y\leq 4$.

The proof of proposition \ref{malena} will be an easy consequence of
the following lemma:
\begin{lemma}
\label{disgiunti}
Let $X$ be Fano and suppose that $\mathcal{G},\mathcal{H}\subset G(\fx)$ 
are two non-empty, disjoint sets of generators 
such that for all $x\in\mathcal{G}$ and 
$y\in\mathcal{H}$, the set $\{x,y\}$ is a primitive collection in
$\fx$.    Then $|\mathcal{G}\cup \mathcal{H}|\leq 4$; moreover, 
if $|\mathcal{G}\cup \mathcal{H}|=4$, $X$ is a toric $S_3$-bundle over
a toric Fano variety, and one of the following two cases occurs:
\begin{enumerate}[(i)]
\item $\mathcal{G}=\{x_1,x_2,x_3\}$, $\mathcal{H}=\{-x_1\}$ with
  primitive relations $(-x_1)+x_1=0$, $(-x_1)+x_2=(-x_3)$ and
  $(-x_1)+x_3=(-x_2)$; 
\item $\mathcal{G}=\{x_1,x_2\}$, $\mathcal{H}=\{-x_1,-x_2\}$ with
  primitive relations $(-x_1)+x_1=0$, $(-x_2)+x_2=0$, $(-x_1)+x_2=v$,
  $x_1+(-x_2)=(-v)$. 
\end{enumerate}
\end{lemma}

\noindent\emph{Remark:} in particular, this implies that the union of $m$
irreducible invariant divisors on $X$ is always a connected set for
$m\geq 5$.  
\begin{dimo}[Proof of lemma \ref{disgiunti}]
We already know from theorem~\ref{structure} that  
$|\mathcal{G}|\leq 3$ and $|\mathcal{H}|\leq 3$: we have to show that if   
$|\mathcal{G}|= 3$, then $|\mathcal{H}|=1 $ (thus by symmetry we also have 
that $|\mathcal{H}|= 3$ implies $|\mathcal{G}|=1$).

We set $\mathcal{G}=\{x_1,x_2,x_3\}$, and let $y\in\mathcal{H}$. 
Since $\{y, x_i\}$ is a primitive collection for $i=1,2,3$, it must be
$y=-x_i$ for some $i$.  
Therefore $\mathcal{H}\subseteq\{-x_1,-x_2,-x_3\}$. Suppose now that
$-x_1\in\mathcal{H}$: then by lemma~\ref{relations}, we have two 
primitive relations $(-x_1)+x_2=(-x_3)$ and $(-x_1)+x_3=(-x_2)$; 
having degree 1, these relations are extremal, hence
the cones $\langle -x_3,x_2\rangle$ and  $\langle -x_2,x_3\rangle$ are
in $\fx$, so $\mathcal{H}=\{-x_1\}$. This gives case $(i)$, and since
$-x_1$ is contained in three primitive collection of order two, 
$X$ is a toric $S_3$-bundle over a toric Fano variety by
theorem~\ref{structure}. 

It remains to show that when $|\mathcal{G}|= |\mathcal{H}|=2 $ the
primitive relations are as described in case $(ii)$. Let
$\mathcal{G}=\{x_1,x_2\}$: if there exists $y\in\mathcal{H}$ such
that $y\neq-x_1$ and $y\neq -x_2$, then we get primitive relations
$x_1+y=(-x_2)$, $x_2+y=(-x_1)$; in particular we see that $y=-x_1-x_2$
is uniquely determined by $x_1$ and $x_2$,
and the cones $\langle x_1,-x_2\rangle$,
$\langle -x_1,x_2\rangle$ are in $\fx$: hence
$\mathcal{H}=\{y\}$. Therefore, the only possibility if
$|\mathcal{H}|=2 $ is that $\mathcal{H}=\{-x_1,-x_2\}$.
Clearly in this case we have primitive relations as in $(ii)$.
Moreover, since $x_1+v=x_2$, also $\{x_1,v\}$ is a primitive
collection: then $x_1$ is contained in three primitive collections of
order 2, and again by theorem~\ref{structure} $X$ is a toric
$S_3$-bundle over a toric Fano variety. 
\end{dimo}
\begin{dimo}[Proof of proposition \ref{malena}]
We remark that $\mathcal{G}_{\eta}$ is non-empty, because 
$f(V(x))=V(\eta)$ means that
$x\in\RelInt\eta$: then
$x$ is in the interior of $\Lambda_{\eta}$, so $x\in\mathcal{G}_{\eta}$. 

Let $y\in\mathcal{G}_{\eta}$ and $z\in\mathcal{H}_{\eta}$: the cone
$\langle y,z\rangle$ crosses the boundary of
$\Lambda_{\eta}$, which is composed of cones in $\fy$. Since $\fx$ is a
subdivision of $\fy$, it must be
$\langle y,z\rangle\not\in\fx$. Hence  $\{y,z\}$ is a
primitive collection for $\fx$. 

Therefore, either $\mathcal{H}_{\eta}=\emptyset$, or
lemma~\ref{disgiunti} applies to $\mathcal{G}_{\eta}$ and
$\mathcal{H}_{\eta}$, and we get (1) and (3). Also (2) is clear, because
if we have a primitive relation $z_1+z_2=v$, then it has degree 1, so
it is extremal, and the cones $\langle z_1,v\rangle$ and $\langle
z_2,v\rangle$ are in $\fx$. Since they can not cross the boundary of
$\Lambda_{\eta}$, it must be $v\in\partial \Lambda_{\eta}$.
\end{dimo}
\begin{dimo}[Proof of corollary \ref{puh}]
Let $p=V(\eta)\in Y$.
Since $\dim\eta=n$, $\fy$ must have at least one generator outside
$\eta$, hence $\mathcal{H}_{\eta}\neq \emptyset$; so
 $|\mathcal{G}_{\eta}\cup\mathcal{H}_{\eta}|\leq 4$ by
 theorem~\ref{malena}.  Notice that $|\mathcal{G}_{\eta}|$ is exactly
 the number of the components of the exceptional divisor that are
 contracted to $p$.

We claim that the cases $|\mathcal{G}_{\eta}|=3$, 
$|\mathcal{H}_{\eta}|=1$ and
$|\mathcal{G}_{\eta}|=|\mathcal{H}_{\eta}|=2$ can not happen. Indeed,
since $\dim\eta=n$, we have that $\Star\eta=\eta$ is strictly
convex. Now, looking carefully at the primitive relations given by $(i)$
and $(ii)$ of lemma~\ref{disgiunti}, one can easily see that 
if $|\mathcal{G}_{\eta}|=3$, $|\mathcal{H}_{\eta}|=1$ or
$|\mathcal{G}_{\eta}|=|\mathcal{H}_{\eta}|=2$, then there exists some
generator $v$ such that $v,-v\in\partial\eta$, which is impossible.

Therefore, we have in any case that $|\mathcal{G}_{\eta}|\leq 2$, 
and if $|\mathcal{G}_{\eta}|=2$ then it must be
 $|\mathcal{H}_{\eta}|=1$, which clearly implies $Y\simeq\pr{n}$.
\end{dimo}


From now on we will consider the case where the morphism 
$f\colon X\rightarrow Y$ is a smooth equivariant blow-up.  
We want to apply the preceeding results 
 to study under which conditions it is possible
that $Y$ is non-projective.
In fact, \emph{as far as we know, there are no known examples of a 
non-projective toric variety $Y$ that becomes Fano after a
smooth equivariant blow-up $X\rightarrow Y$}. 

Let $A\subset Y$ be the center of the blow-up. 
If $A$ is a point, $Y$ is always projective; if $\dim A=n-2$, a line in a
non-trivial fiber of the blow-up has anticanonical degree 1, so it
is extremal by proposition~\ref{gradouno}, and again $Y$ is projective.

The result we obtain is the following:
\begin{teo}
\label{sushi}
Let $Y$ be a smooth complete non-projective toric variety of dimension
$n=\dim Y\geq 4$, and $A\subset Y$ an irreducible invariant
subvariety. If the
blow-up of $Y$ along $A$ is Fano, then $3\leq\dim A \leq n-3$, $n\geq
6$, $\rho_Y-\rho_A\leq 2$ and $\rho_A >1$.
\end{teo}
\noindent Theorem~\ref{sushi} will be a consequence of
proposition~\ref{picard} and proposition~\ref{ciop}.


\begin{prop}
\label{picard}
Let $X$ be an $n$-dimensional toric Fano variety, $f\colon
X\rightarrow Y$ a smooth equivariant blow-up and  $E\subset X$ 
the exceptional divisor of $f$. If $\rho_X-\rho_E=3$, or if
$\rho_X-\rho_E=2$ and $X$ is a toric $S_2$-bundle $\pi\colon X\rightarrow Z$ 
with $\pi_{|E}\colon E\rightarrow Z$ a $\pr{1}$-bundle,
then $Y$ is projective.
\end{prop}
\begin{dimo}
We will show that $Y$ is projective when $\rho_X-\rho_E=3$; 
in the same way one can prove the other case. 
By \ref{structure}.1, there is a toric $S_3$-bundle $\pi\colon
X\rightarrow Z$ over a toric Fano
variety $Z$. Let $E=V(x)$: keeping the notation of theorem~\ref{structure}, 
we know that
in $\fx$ there are the nine primitive relations:
\[
\begin{array}{c} 
x+(-x)=0\ \ y+(-y)=0\ \ w+(-w)=0\ \ x+y=(-w)\ \ x+w=(-y)\\
(-x)+(-y)=w\ \ y+w=(-x)\ \ (-x)+(-w)=y\ \ (-y)+(-w)=x. 
\end{array}
\]
Let $\mathcal{K}\subset \PC(\fx)$ be the set of these nine primitive
collections; we also know that there is a bijection between
$\PC(\fx)\smallsetminus\mathcal{K}$ and $\PC(\f_Z)$, induced by the
surjective morphism $\pi_*\colon\NE(X)\rightarrow\NE(Z)$. 

We claim that for every primitive collection
$P\in\PC(\fx)\smallsetminus\mathcal{K}$, $r(P)$ is extremal in
$\NE(X)$ if and only if $\pi_*(r(P))$ is extremal in $\NE(Z)$. Indeed,
clearly if $r(P)$ decomposes as a sum of other primitive collections,
the decomposition holds for $\pi_*(r(P))$, and conversely.

Now consider the class $\omega=r(Q)$ coming from the blow-up $f$; 
$Y$ is projective if and only if $\omega$ is extremal in
$\NE(X)$. If $\omega$ is in 
$\mathcal{K}$, then it has degree 1, hence it is extremal.
So we can suppose that $\omega\in\PC(\fx)\smallsetminus\mathcal{K}$;
we want to show that $\pi_*(\omega)$ is extremal in $\NE(Z)$.
We remark that $\omega_{|E}$ is extremal in $E$, because it
corresponds to a $\pr{r}$-bundle on the center of the blow-up
$f$. But this implies that $\pi_*(\omega)$ is extremal in $Z$, because
it is the image of $\omega_{|E}$ under the restriction $\pi_{|E}\colon
E\rightarrow Z$, which is a $\pr{1}$-bundle. 
\end{dimo}

In the case $\dim A\leq 2$ it is sufficient to suppose that
$Y$ is not Fano to get strong conditions on $A$. This result is due
to J.~Wi\'sniewski~\cite{wisn}, and it is particularly
interesting because it holds in a general (non-toric) setting:
\begin{teo}[\cite{wisn}, 3.5 and 3.6]
\label{wis}
Let $X$ be a Fano manifold of dimension $n\geq 4$ and $f\colon
X\rightarrow Y$ the blow-up of a smooth variety $Y$ along a smooth
subvariety $A$, with $\dim A\leq 2$; let $E$ be the exceptional
divisor of $f$. 
Suppose that $Y$ is not Fano. Then one of the following cases occurs:
\begin{enumerate}[(i)]
\item
 $A\simeq\pr{1}$,
$N_{A/Y}\simeq\mathcal{O}(-1)^{\oplus(n-1)}$ and consequently
$E\simeq\pr{1}\times\pr{n-2}$;
\item $A\simeq\pr{2}$ and $N_{A/Y}$ is either
  $\mathcal{O}(-1)^{\oplus(n-2)}$ and $n\geq 5$, or 
$\mathcal{O}(-2)\oplus\mathcal{O}(-1)^{\oplus(n-3)}$,
or $T\pr{2}(-3)\oplus\mathcal{O}(-1)^{\oplus(n-4)}$, where
$T\pr{2}$ is the tangent bundle to $\pr{2}$;
\item $A\simeq\pr{1}\times\pr{1}$ and
 $N_{A/Y}\simeq\mathcal{O}(-1,-1)^{\oplus(n-2)}$ and consequently
$E\simeq\pr{1}\times\pr{1}\times\pr{n-3}$;
\item $A$ is a surface with a
 $\pr{1}$-bundle structure $A\rightarrow C$ over a
  smooth curve $C$. In this case $X$ admits another blow-down
  $X\rightarrow Y_1$  which contracts the divisor $E$ to a smooth
  codimension 2 subvariety of $Y_1$.
\end{enumerate}
\end{teo}
\begin{prop}
\label{ciop}
Let  $X\rightarrow Y$ be the blow-up of a smooth complete non-Fano
toric variety  $Y$ along a smooth invariant surface $A\subset Y$; let
$E\subset X$ be the exceptional divisor.
 Suppose that $X$ is Fano, $\rho_X=5$ and $\rho_A=2$. 
Then $Y$ is projective. 

Moreover, 
$A\simeq\pr{1}\times\pr{1}$ or $A\simeq S_1$, $E\simeq
A\times\pr{n-3}$, and $X$ is either a toric $S_2$-bundle over
$\pr{1}\times\pr{n-3}$, or the blow-up of a toric $\pr{1}$-bundle over
$A\times\pr{n-3}$ along a codimension 2 smooth invariant 
subvariety contained in a
section.
\end{prop}
\begin{dimo}
We know, by theorem~\ref{wis}, that  $A\simeq\mathbb{F}_a$ with
$a\in\mathbb{Z}_{\geq 0}$. We treat the case $a>0$; the case $a=0$ is
very similar and easier, because the normal bundle $N_{E/X}$ is known.

Let $E=V(x)$: 
we know by theorem~\ref{wis} that $E$ has a $\pr{1}$-bundle
structure over an $(n-2)$-dimensional toric variety. Hence in $\f_E$ we
have primitive relations:
\[ u_0+u_1=0,\quad v_0+v_1=au_0,\quad   z_0+\cdots+z_{n-3}=0, \]
while in $\fx$ we know the relations 
$u_0+u_1=x$ and  $z_0+\cdots+z_{n-3}=x$. Moreover, either
$\{v_0,v_1\}$ is a primitive collection in $\fx$, or $\{x,v_0,v_1\}$
is.  Since $-K_E$ is nef, it must be $a=1$ or $a=2$.

We have $\rho_E=3$ and $\rho_X=5$, thus $x$ is contained in
exactly two primitive collections of order two $\{x,w_0\}$ and
$\{x,w_1\}$, and 
\[ G(\fx)=\{x,u_0,u_1,v_0,v_1,w_0,w_1,z_0,\ld,z_{n-3}\}. \]
We examine the possible primitive relations of $\{x,w_0\}$ and
$\{x,w_1\}$.

\medskip

\noindent 1) $x+w_0=0,\quad x+w_1=z_0$. 

\vspace{15pt} 

\hspace{150pt}\psfig{figure=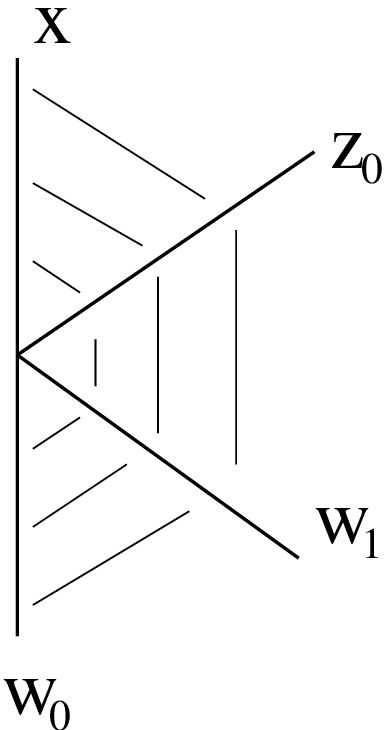,width=1.5cm}
\vspace{15pt}

We remark that if $\{x,v_0,v_1\}$ is a
 primitive collection in $\fx$, then by proposition~\ref{idea} the
 associated relation should be $x+v_0+v_1=p_0+p_1$. Since the degree
 is 1, we would get $\langle x,p_0,p_1\rangle\in\fx$, hence this would
 give in $\f_E$ the primitive relation $v_0+v_1=p_0+p_1$, a contradiction.
 Hence  $\{v_0,v_1\}$ is a primitive collection in $\fx$, with
 relation $v_0+v_1=au_0$, so $a=1$.
Then we get primitive relations:
\[ z_0+w_0=w_1,\quad  w_1+z_1+\cdots+z_{n-3}=0,\quad v_0+v_1=u_0.
\]
We easily see that $X$ is the blow-up along $V(\langle z_0,w_0\rangle)$ of the
variety $X'$ given by the primitive relations:
\[ x+w_0=0,\quad v_0+v_1=u_0,\quad u_0+u_1=x,\quad z_0+\cdots+z_{n-3}=x. \]
$X'$ is a $\pr{1}$-bundle over $\pr{n-3}\times S_1$.
Clearly the relation $z_0+\cdots+z_{n-3}=x$ is extremal in $\NE(X)$, thus
$Y$ is projective.

\medskip

\noindent 2) $x+w_0=0,\quad x+w_1=v_1$. Exactly as in the preceeding case
we see that $a=1$, $v_0+v_1=u_0$ is a primitive relation in $\fx$ and
$X$ is the blow-up of $X'$ along $V(\langle v_1,w_0\rangle)$. 

\medskip

\noindent 3) $x+w_0=0,\quad x+w_1=u_1$. 
As in case 1), we have that
$a=1$ and $v_0+v_1=u_0$ is a primitive relation in $\fx$. Moreover we know the
primitive relations $u_0+w_1=0$ and $ u_1+w_0=w_1$. We are in the
hypothesis of \ref{structure}.2a, 
and $x$ can not appear in other
primitive collections, because this would give other primitive
collections in $E$: hence $X$ is toric $S_2$-bundle over
$\pr{n-3}\times\pr{1}$. In particular, $Y$ is projective, by
proposition~\ref{picard}.

\medskip

\noindent 4) $x+w_0=0,\quad x+w_1=u_0$. This is exactly as  the
preceeding case, except for the fact that the primitive relation
associated to $\{v_0,v_1\}$ can be either $v_0+v_1=u_0$ or 
$v_0+v_1=w_1$. In both cases $X$ is a toric $S_2$-bundle over
$\pr{n-3}\times\pr{1}$ and $Y$ is projective.

\medskip

\noindent 5) $x+w_0=u_0,\quad x+w_1=u_1$. 

\vspace{15pt} 

\hspace{110pt}\psfig{figure=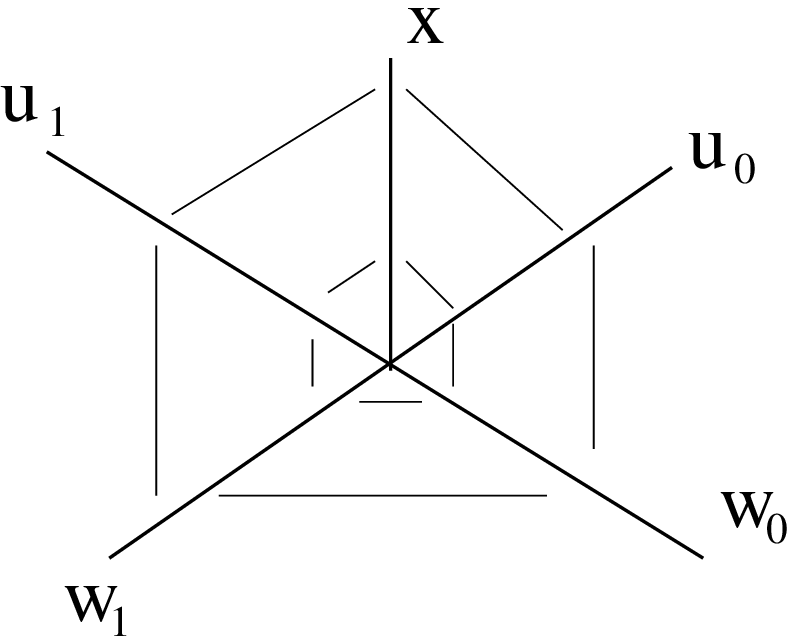,width=3.5cm}
\vspace{15pt}

This is the only case with
$-x\not\in G(\fx)$, and we are again in the hypothesis of
\ref{structure}.2a. 
This implies that  $E$ is Fano, hence $a=1$. 
If $\{x,v_0,v_1\}$ is a primitive collection in $\fx$, the 
primitive relation could be $x+v_0+v_1=u_0$ or
 $x+v_0+v_1=w_0$, but both cases are impossible: indeed, we have seen in the
 proof of \ref{structure}.2a 
that any primitive relation involving $x$
 with coefficient 1 can not contain any other generator of the fan of
 $S_2$. 
Therefore we get the same as case 4), $X$ is toric $S_2$-bundle over
$\pr{n-3}\times\pr{1}$ and $Y$ is projective.
\end{dimo}
\begin{dimo}[Proof of theorem \ref{sushi}]
Since $Y$ is non-projective, it must be $\rho_Y\geq 4$ (see
\cite{kleinschsturm}). Moreover, by proposition \ref{picard} we have 
$\rho_Y-\rho_A\leq 2$, hence $\rho_A\geq 2$.
So it clearly can not be $\dim A=1$. If $\dim A=2$, the possible cases
are described by theorem~\ref{wis}. Again, it can not be 
$A\simeq\pr{2}$; in the other cases we have 
$\rho_A=2$ and $\rho_X=5$: then
proposition~\ref{ciop} applies and $Y$ is projective.
Hence,  it must be $\dim A\geq 3$.
\end{dimo}

\stepcounter{subsubsec}
\subsubsection{Equivariant birational morphisms in dimension 4:
  analysis of the possible subdivisions} 
\label{quinta} In \cite{sato}, H.~Sato shows that every equivariant
birational 
morphism between toric Fano 3-folds is a composite of smooth
equivariant blow-ups between toric Fano 3-folds. 
We recall that in general an equivariant birational morphism $f$ between
smooth complete toric varieties doesn't admit a factorization in a
sequence of smooth blow-ups, if the dimension is greater than two
($f$ will have a factorization as a sequence of smooth equivariant
blow-ups and blow-downs; see \cite{morelli} and  \cite{abr1}).
In case of an equivariant birational morphism between toric Fano
3-folds, not only is there a factorization in smooth equivariant
blow-ups, but also all the intermediate 3-folds are Fano.
Such a result is proven in two steps:

\noindent 1) if the source is Fano, then the morphism is a composite
of smooth equivariant blow-ups;

\noindent 2) if also the target is Fano, then all the intermediate
3-folds are Fano.

\noindent In \cite{cras} the author has shown, with an explicit 
counterexample, that the same result does not hold in dimension 4:
\begin{prop}[\cite{cras}, 3.1]
There exist two toric Fano 4-folds $X$ and $Y$ and a birational
equivariant morphism $f\colon X\rightarrow Y$, such that $f$ doesn't
admit a decomposition in smooth equivariant blow-ups between toric
Fano 4-folds.
\end{prop}
\noindent In this section we show that at least 1) holds in
dimension 4:
\begin{teo}
\label{facto}
Let $X$, $Y$ be two toric 4-folds and $f\colon X\rightarrow Y$ an
equivariant birational morphism. Suppose that $X$ is Fano. Then $f$
factorizes as a sequence of smooth equivariant blow-ups between toric
4-folds.
\end{teo}
\noindent The 
proof of theorem \ref{facto} is analogous to Sato's proof in
dimension 3: we fix a maximal cone $\sigma\in \Sigma_Y$ and
 study the possible subdivisions of $\sigma$ in $\fx$, applying the results 
of the preceeding section about the possible positions of new generators.
Geometrically, this means that 
we are considering the restriction of
$f$ to $f^{-1}(\mathcal{U}_{\sigma})$, where
$\mathcal{U}_{\sigma}\subset Y$ is an invariant open subset isomorphic
to $\mathbb{A}^n$.
The result we obtain is the following:
\begin{prop}
\label{subdiv}
Let $X$ and $Y$ be two toric 4-folds and $f\colon X\rightarrow Y$ an
equivariant birational morphism. Suppose that $X$ is Fano and let
$\sigma=\langle y_1,y_2,y_3,y_4\rangle\in\fy$. Then the possible
subdivisions of $\sigma$ in $\fx$ are 17 and are all given by a
sequence of at most 3 star-subdivisions. These 17 subdivisions are
described in the following list, where we give for each subdivision: 
the sequence of cones whose
star-subdivisions (in the given order) give the subdivision;  the
primitive relations inside the subdivision.
The last four cases can occur only if $Y\simeq\pr{4}$.    

\medskip

{\upshape \noindent 1) $\sigma\in\fx$, i.~e.\ $\sigma$ is not subdivided.\vspace{4pt}\\
\noindent 2) 
$\langle y_1,y_2\rangle$; the only primitive relation is
  $y_1+y_2=x_1$.\vspace{4pt}\\
\noindent 3)
$\langle y_1,y_2,y_3\rangle$; the only primitive relation is
  $y_1+y_2+y_3=x_1$.\vspace{4pt}\\
\noindent4)
 $\langle y_1,y_2,y_3,y_4\rangle$; the only primitive relation is
  $y_1+y_2+y_3+y_4=x_1$.\vspace{4pt}\\
\noindent 5)
 $\langle y_1,y_2\rangle$, $\langle y_3,y_4\rangle$; 
  primitive relations: 
  $y_1+y_2=x_1$,  $y_3+y_4=x_2$.\vspace{4pt}\\
\noindent 6)
 $\langle y_1,y_2,y_3\rangle$, $\langle y_1,x_1\rangle$; 
  primitive relations: 
  $y_1+y_2+y_3=x_1$,  $y_1+x_1=x_2$, $y_2+y_3+x_2=2x_1$.\vspace{4pt}\\
7)
 $\langle y_1,y_2,y_3\rangle$, $\langle y_1,y_4\rangle$; 
  primitive relations:
  $y_1+y_2+y_3=x_1$,  $y_1+y_4=x_2$, $y_2+y_3+x_2=x_1+y_4$.\vspace{4pt}\\
\noindent 8)
 $\langle y_1,y_2,y_3\rangle$, $\langle y_1,y_2\rangle$; 
  primitive relations:
  $y_1+y_2=x_2$,  $x_2+y_3=x_1$.\vspace{4pt}\\
9)
 $\langle y_1,y_2,y_3,y_4\rangle$, $\langle y_1,y_2,y_3\rangle$; 
  primitive relations:
  $y_1+y_2+y_3=x_2$,  $x_2+y_4=x_1$.\vspace{4pt}\\
10)
 $\langle y_1,y_2,y_3\rangle$, $\langle y_2,y_3,y_4\rangle$,
  $\langle y_4,x_2\rangle$;
  primitive relations:
  $y_1+y_2+y_3=x_2$,  $y_2+y_3+y_4=x_3$, $y_1+x_3=x_1$, $y_4+x_2=x_1$,
  $y_2+y_3+x_1=x_2+x_3$.\vspace{4pt} \\
11)
 $\langle y_1,y_2,y_3,y_4\rangle$, $\langle y_1,y_2\rangle$; 
 primitive relations: 
  $y_1+y_2=x_2$,  $y_3+y_4+x_2=x_1$.\vspace{4pt}\\
12)
 $\langle y_1,y_2,y_3,y_4\rangle$, $\langle y_1,y_2\rangle$,
  $\langle y_3,y_4\rangle$;
  primitive relations: 
  $y_1+y_2=x_2$,  $y_3+y_4=x_3$, $x_2+x_3=x_1$.\vspace{4pt}\\ 
13)
$\langle y_1,y_2,y_3,y_4\rangle$, $\langle y_1,y_2\rangle$,
  $\langle x_2,y_3\rangle$;
  primitive relations:
  $y_1+y_2=x_2$,  $y_3+x_2=x_3$, $y_4+x_3=x_1$.\vspace{4pt}\\
14)
 $\langle y_1,y_2,y_3,y_4\rangle$, $\langle y_1,x_1\rangle$; 
  primitive relations:
  $y_1+y_2+y_3+y_4=x_1$,  $y_1+x_1=x_2$, $x_2+y_2+y_3+y_4=2x_1$.\vspace{4pt}\\
15)
 $\langle y_1,y_2,y_3,y_4\rangle$, $\langle y_1,y_2,y_3\rangle$,
  $\langle y_1,x_1\rangle$; 
  primitive relations:
  $y_1+y_2+y_3=x_2$,  $x_2+y_4=x_1$, $y_1+x_1=x_3$, $y_2+y_3+x_3=x_1+x_2$.\vspace{4pt}\\
16)
 $\langle y_1,y_2,y_3,y_4\rangle$, $\langle y_1,y_2\rangle$,
  $\langle y_3,x_1\rangle$;
  primitive relations:
  $y_1+y_2=x_2$,  $y_3+y_4+x_2=x_1$,  $y_3+x_1=x_3$, $y_4+x_2+x_3=2x_1$.\vspace{4pt}\\
17)
$\langle y_1,y_2,y_3,y_4\rangle$, $\langle y_1,y_2\rangle$,
  $\langle x_1,x_2\rangle$;
  primitive relations:
  $y_1+y_2=x_2$,  $y_3+y_4+x_2=x_1$, $x_1+x_2=x_3$, $y_3+y_4+x_3=2x_1$.}
\vspace{15pt}

\hspace{5pt}(1)\psfig{figure=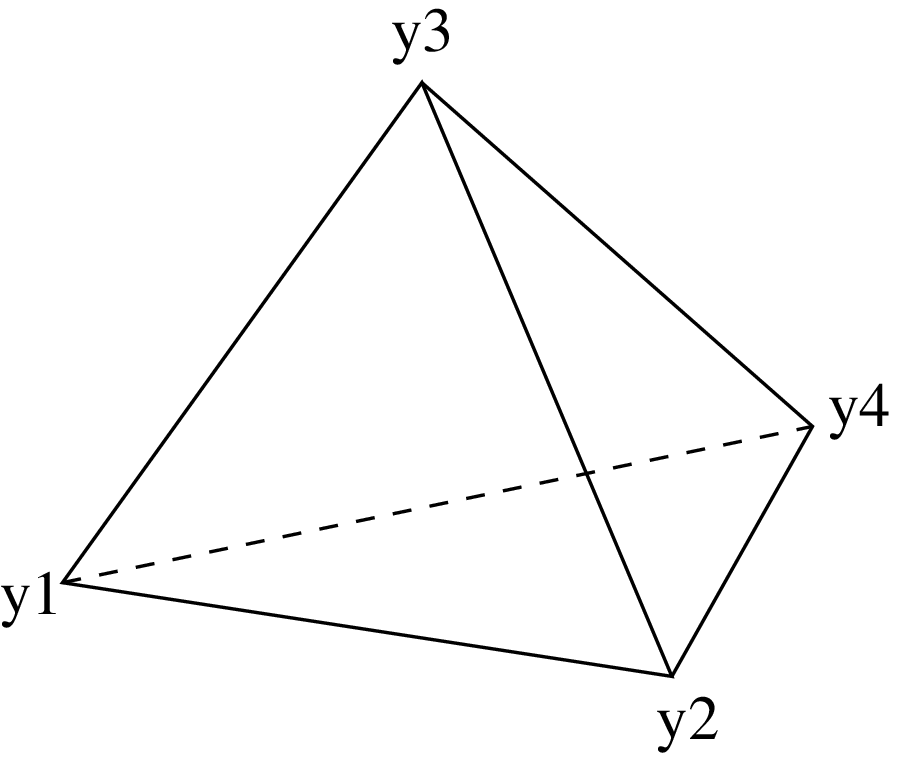,width=5cm}\hspace{15pt}(2)\psfig{figure=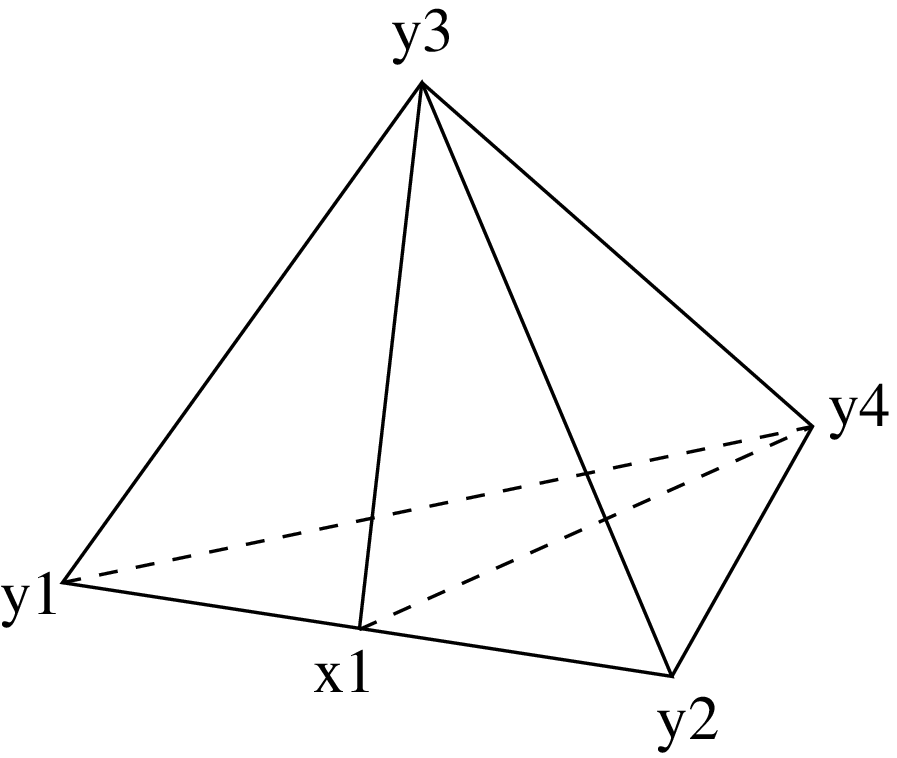,width=5cm}
\vspace{15pt}

\hspace{5pt}(3)\psfig{figure=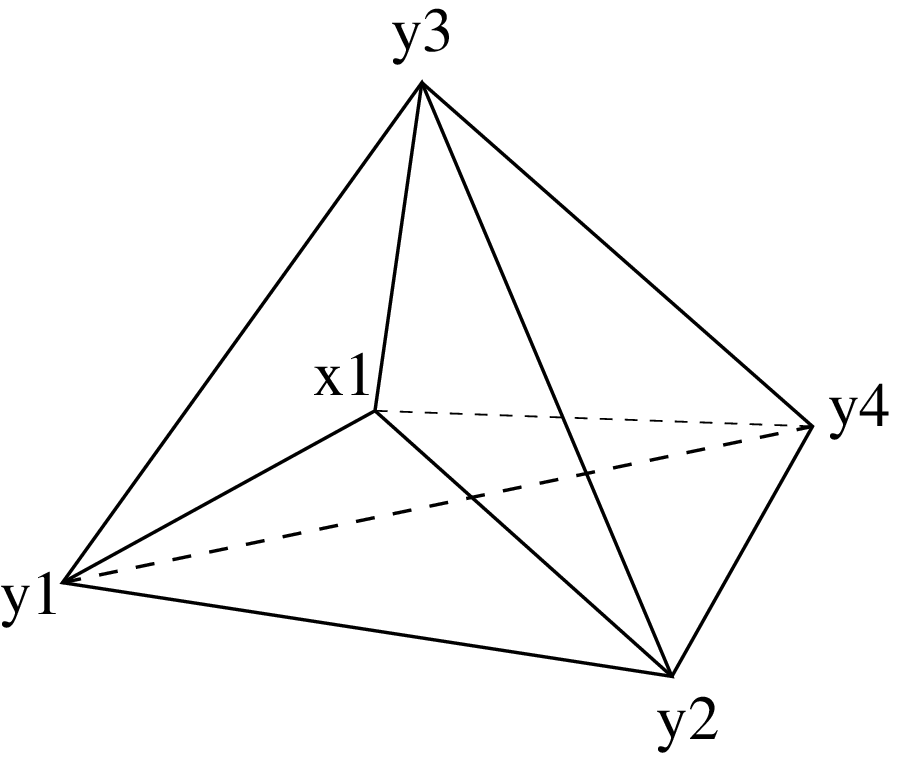,width=5cm}\hspace{15pt}(4)\psfig{figure=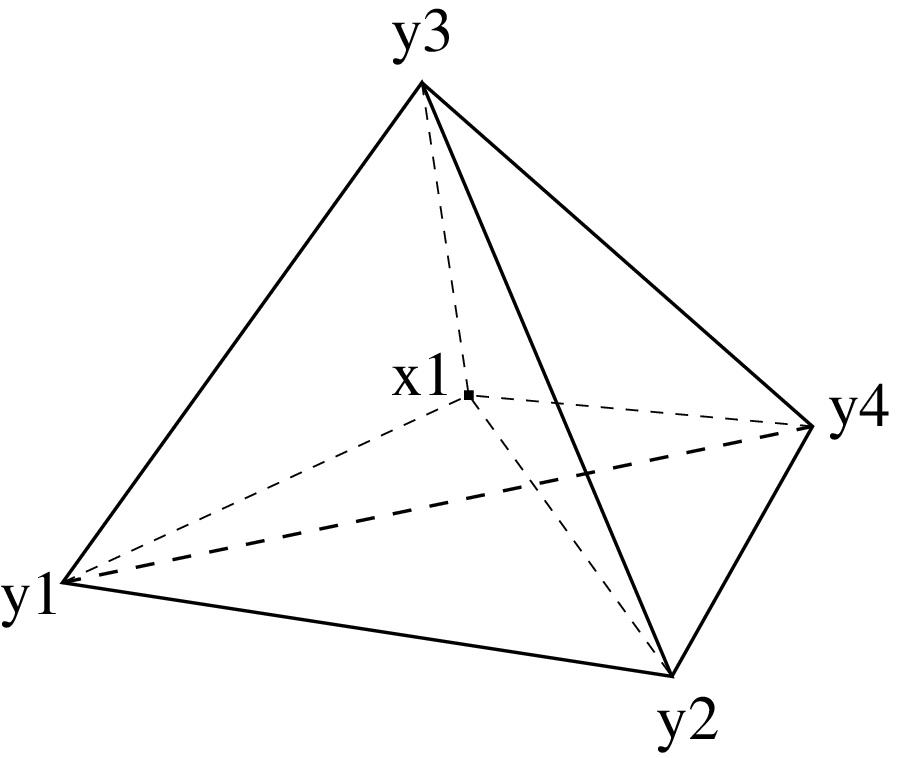,width=5cm}

\vspace{15pt}

\hspace{5pt}(5)\psfig{figure=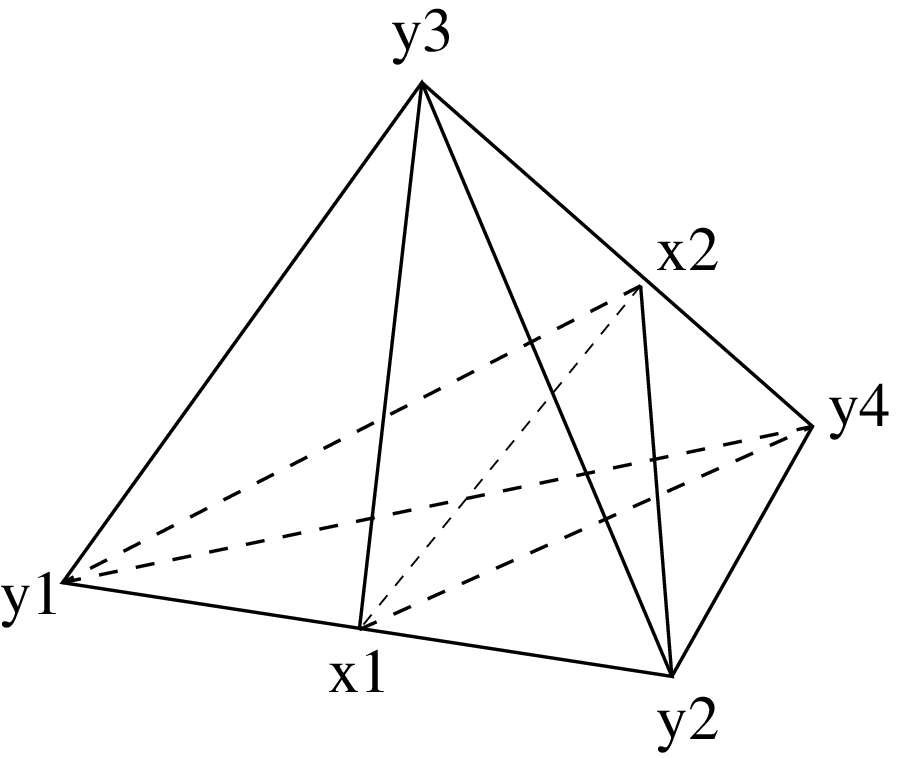,width=5cm}\hspace{15pt}(6)\psfig{figure=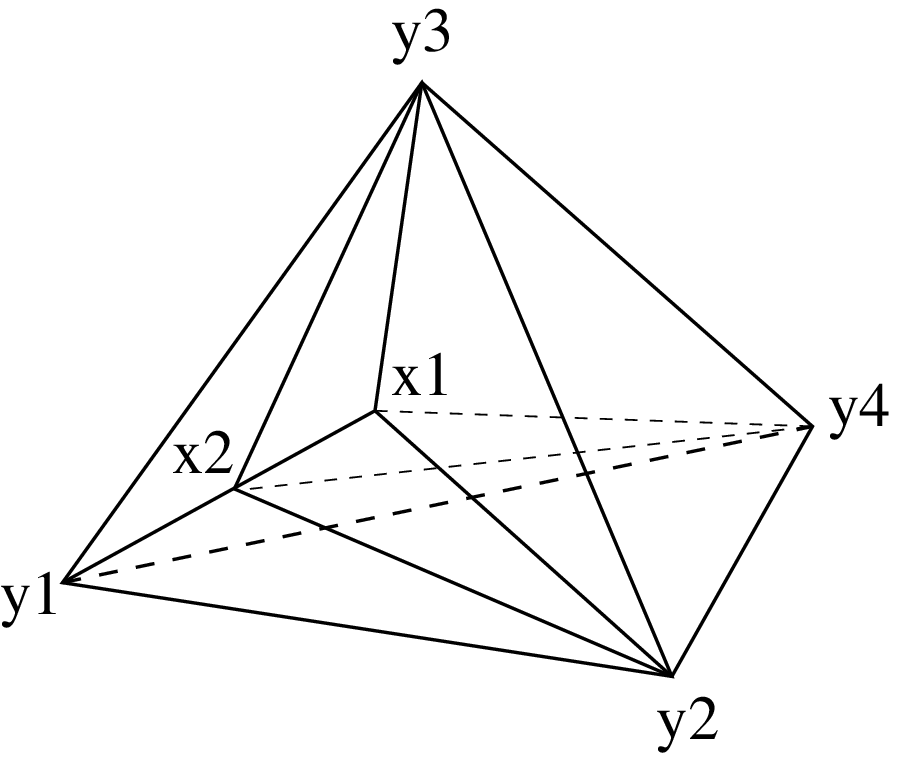,width=5cm}

\vspace{15pt} 

\hspace{5pt}(7)\psfig{figure=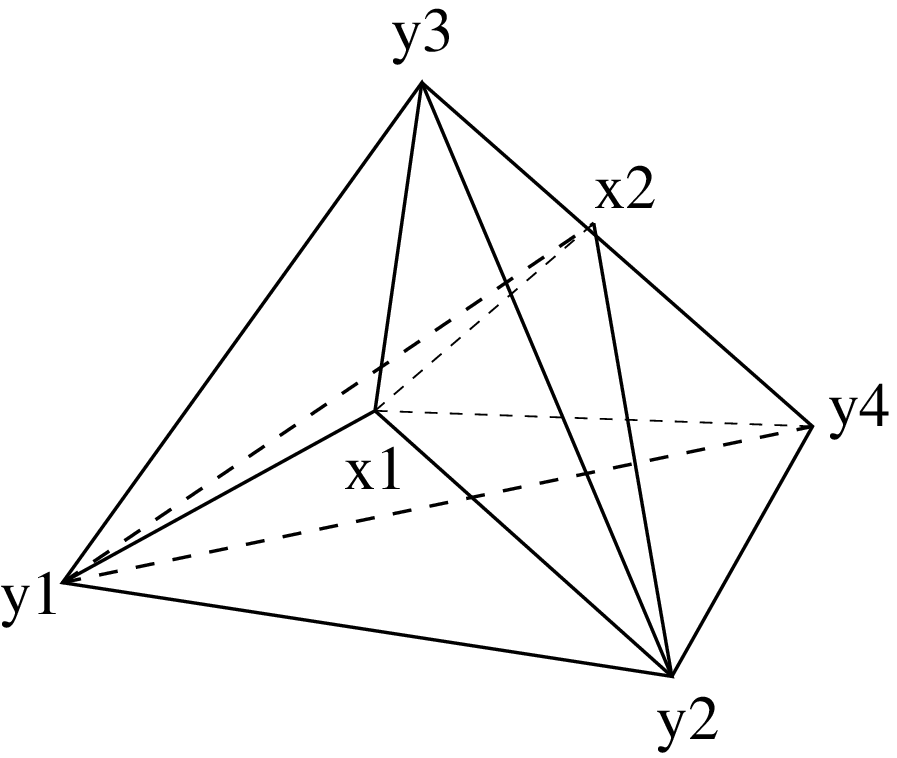,width=5cm}\hspace{15pt}(8)\psfig{figure=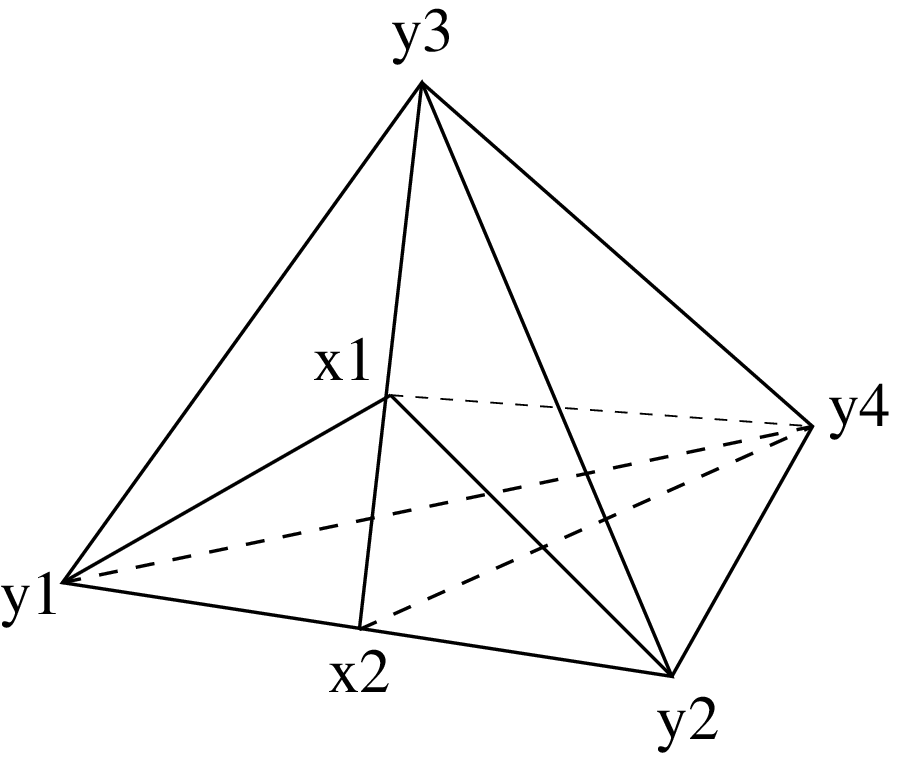,width=5cm}

\vspace{15pt}

\hspace{5pt}(9)\psfig{figure=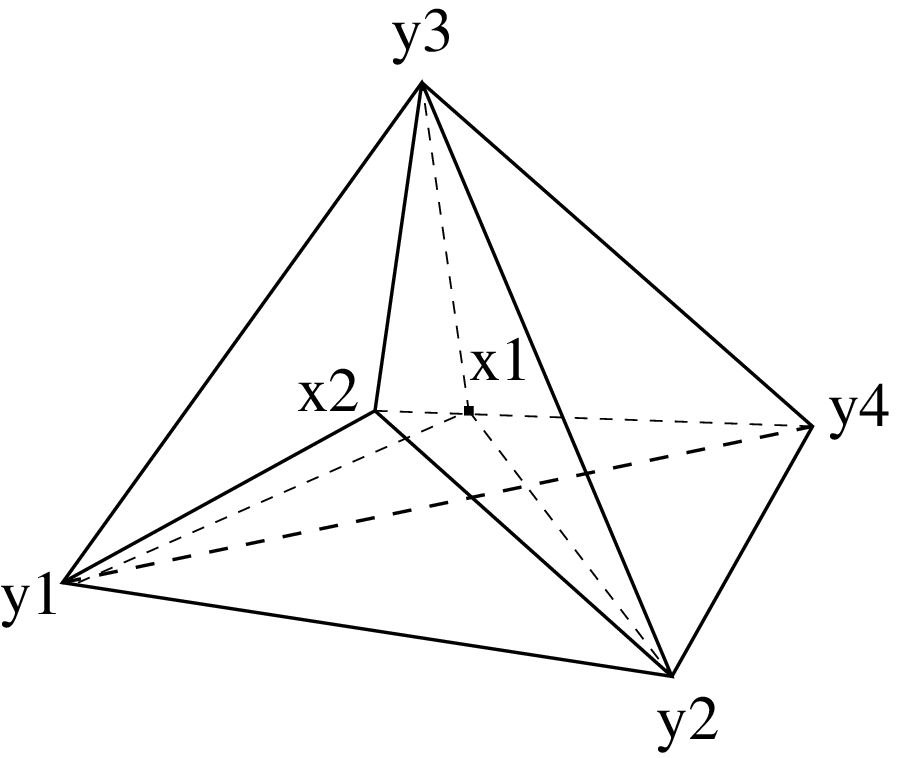,width=5cm}\hspace{15pt}(10)\psfig{figure=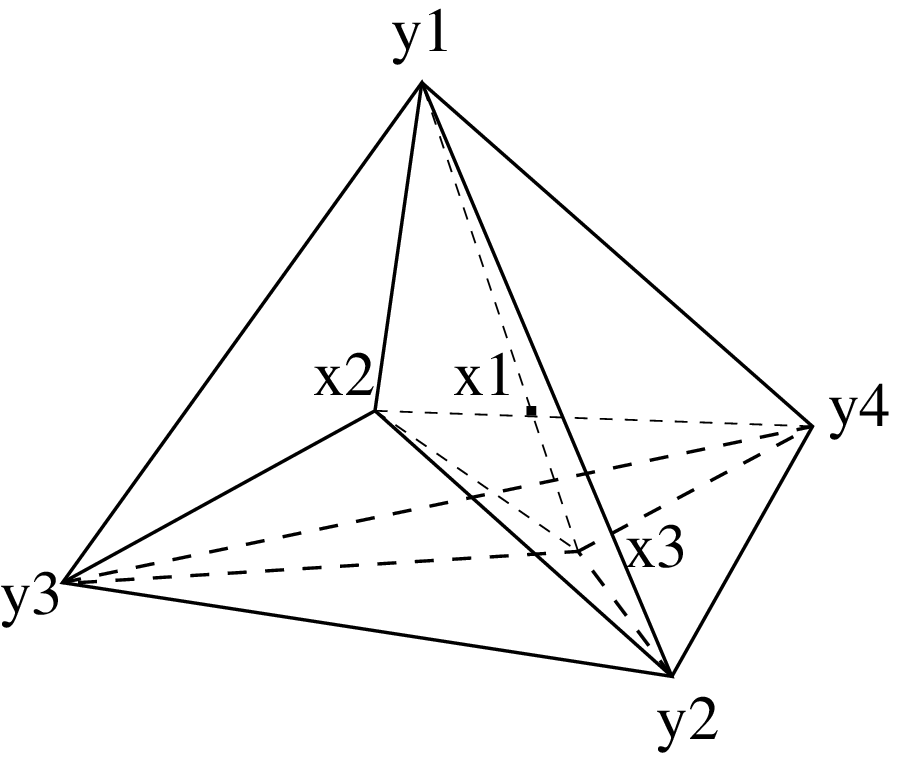,width=5cm}

\vspace{15pt} 

\hspace{5pt}(11)\psfig{figure=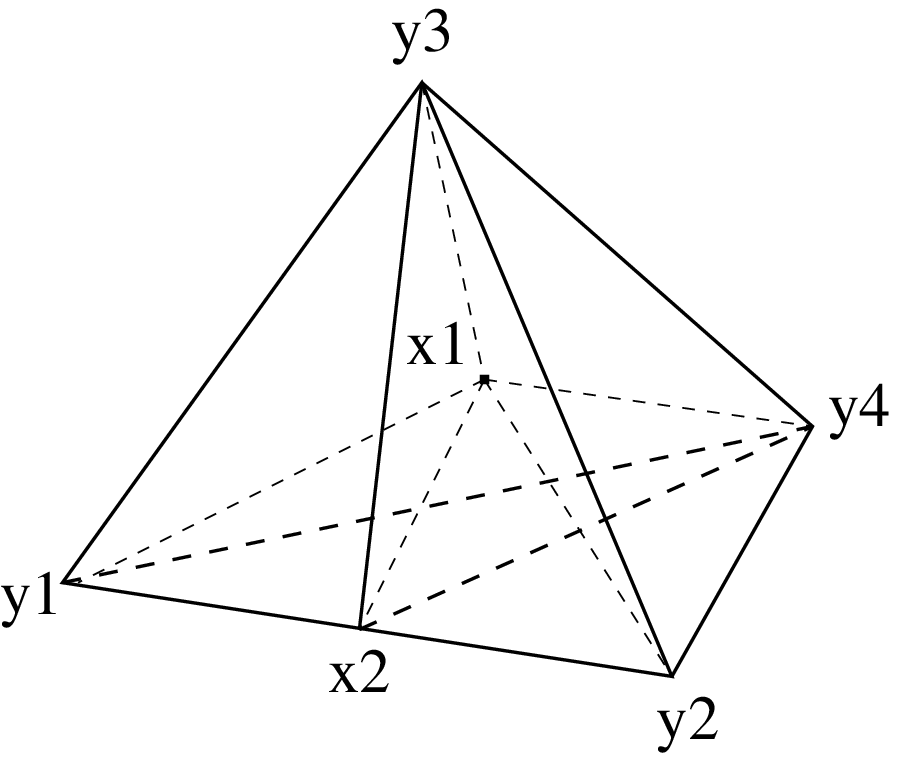,width=5cm}\hspace{15pt}(12)\psfig{figure=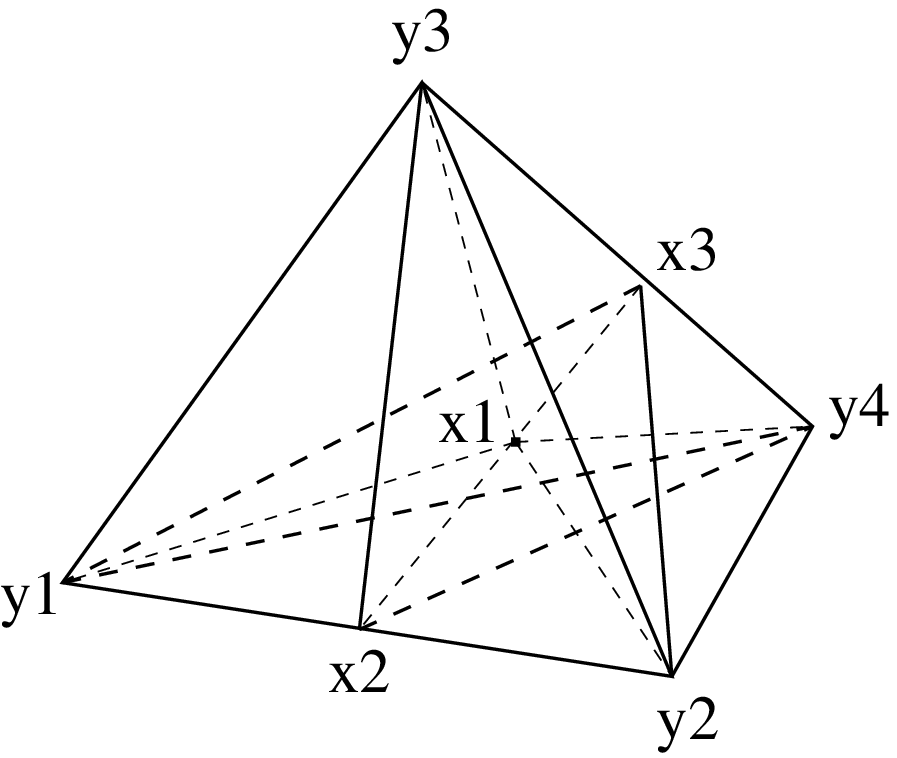,width=5cm}

\vspace{15pt}

\hspace{5pt}(13)\psfig{figure=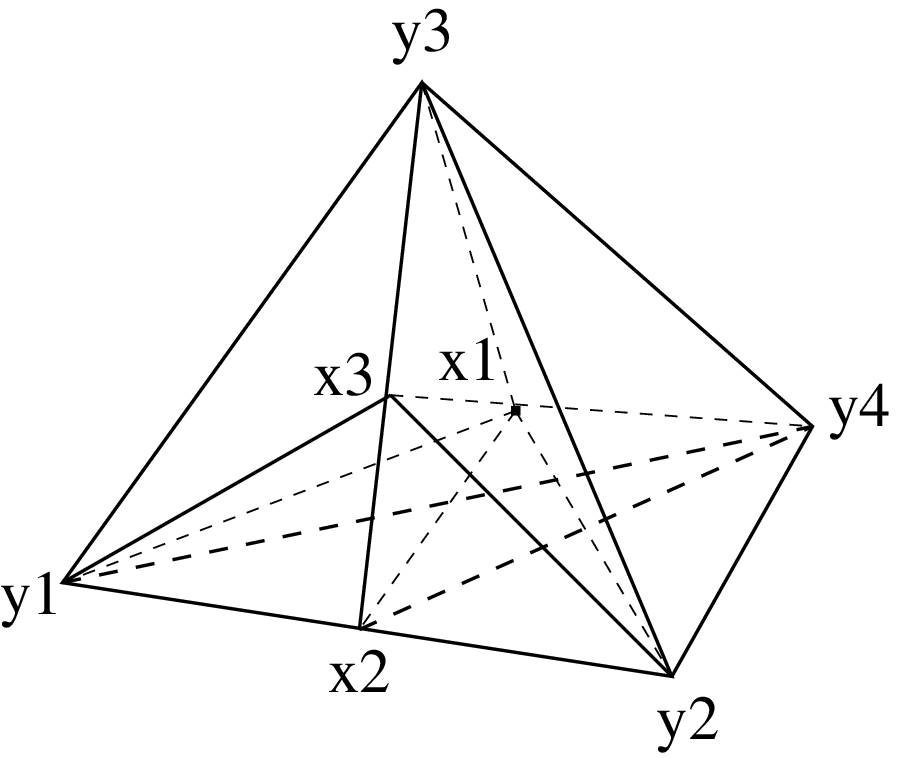,width=5cm}\hspace{15pt}(14)\psfig{figure=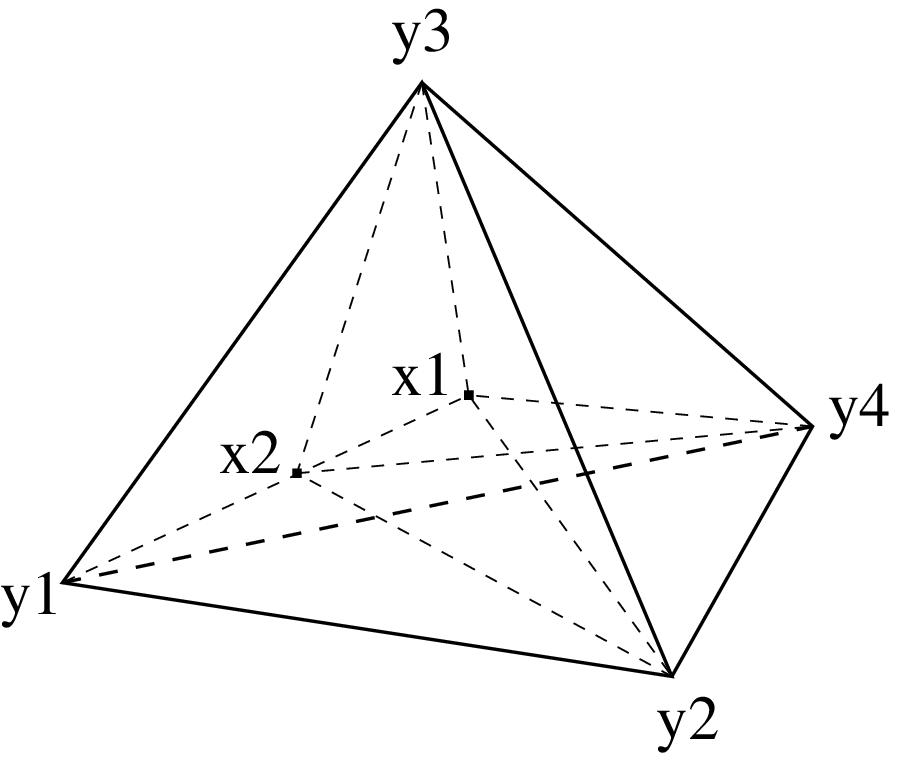,width=5cm}

\vspace{15pt} 

\hspace{5pt}(15)\psfig{figure=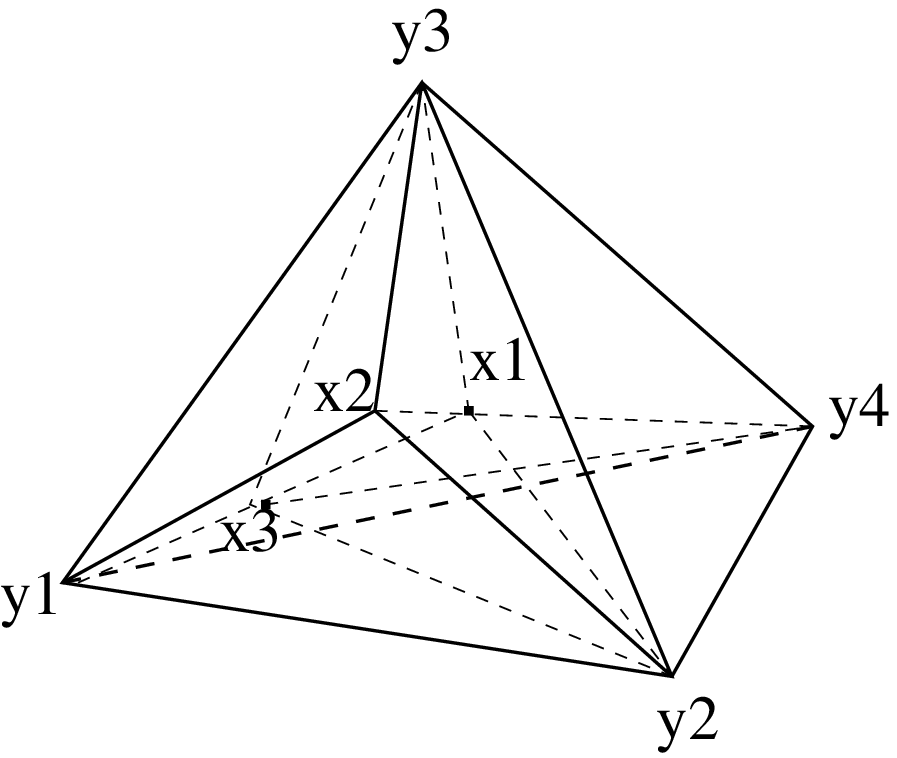,width=5cm}\hspace{15pt}(16)\psfig{figure=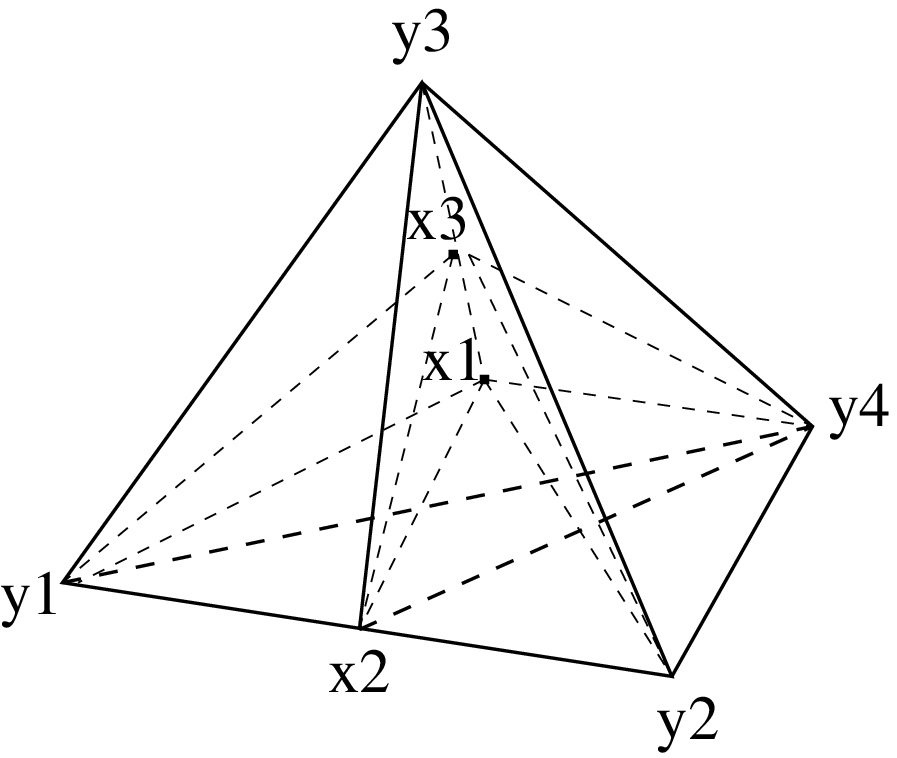,width=5cm}

\vspace{15pt}

\hspace{30pt}(17)\psfig{figure=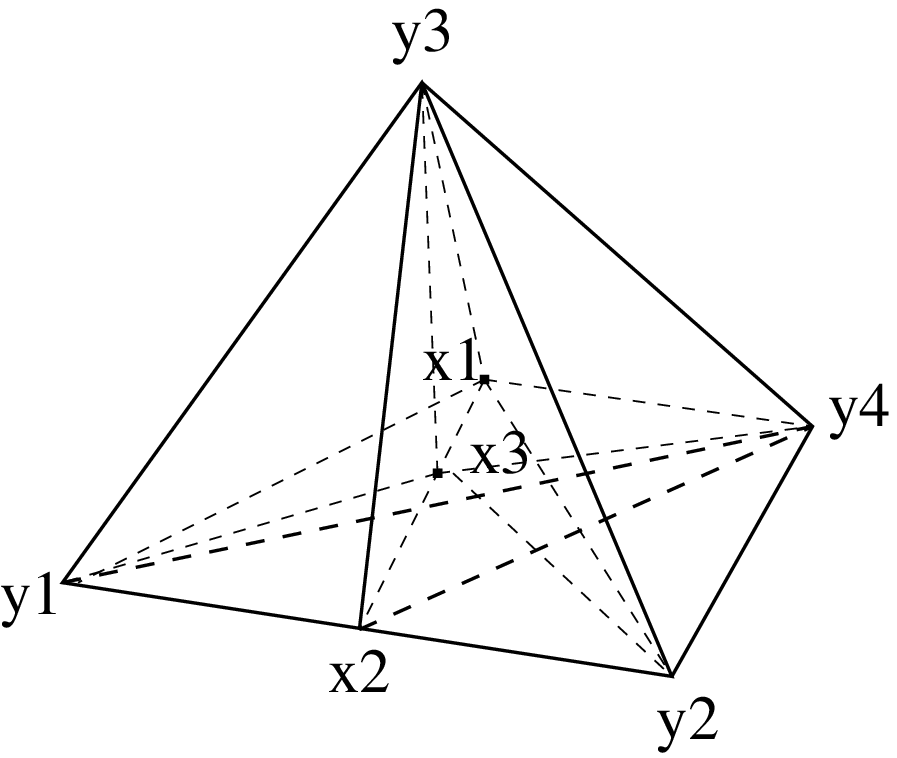,width=5cm}
\end{prop}

\medskip

\noindent\emph{Remark:}  in the case $Y\simeq\pr{4}$, 
it is easy to see that all these subdivisions can actually occur.

\medskip

The proofs of proposition \ref{subdiv} and theorem \ref{facto} will
take all the rest of the section. 

This first lemma is lemma 5.4 in   \cite{sato}; it is a 
result of linear algebra, we give the proof for completeness:
\begin{lemma}
\label{centro}
Let $X$, $Y$ be two complete, smooth, n-dimensional toric varieties,
and $f\colon X\rightarrow Y$ an equivariant birational morphism. Fix a
cone $\tau=\langle y_1,\ldo,y_m\rangle$ in $\fy$ and let
$\widetilde{\tau}=\langle x_1,\ldo, x_k\rangle$ 
be the unique cone in $\fx$ such that $y_1+\cdots +y_m\in 
\RelInt(\widetilde{\tau})$. Then
there exist a partition $J_1,\ldo,J_k$ of the set $\{1,\ldo,m\}$ such that 
\[
x_i = \sum_{j\in J_i} y_j \qquad \forall i=1,\ldo,k.
\]
\end{lemma}
\begin{dimo}
Since $X$ is complete, the point
$y_1+\cdots +y_m$ actually lies in some cone of $\fx$. Moreover, $\fx$ is
a subdivision of $\fy$, thus it must be
$\widetilde{\tau}\subseteq\tau$. Then for all $i=1,\ldo,k$ we have
\[
x_i=a_{i1}y_1+\cdots+a_{im}y_m
\]
with $a_{ij}\in\mathbb{Z}_{\geq 0}$. On the other hand,
since $y_1+\cdots +y_m \in \RelInt\widetilde{\tau}$, there exist 
$b_1,\ldo,b_k\in\mathbb{Z}_{>0}$ such that
\[
y_1+\cdots +y_m=\sum_{i=1}^k b_i x_i=\sum_{i=1}^k b_i
\left(\sum_{h=1}^m a_{ih} y_h\right)
=\sum_{h=1}^m\left(\sum_{i=1}^k b_i a_{ih}\right)y_h.
\]
Since $Y$ is smooth, 
the set $\{y_1,\ldo,y_m\}$ is a part of a 
basis of the lattice: hence we have
\[
\sum_{i=1}^k b_i a_{ih}=1 \quad \forall h=1,\ldo,m.
\]
Since $b_i\in\mathbb{Z}_{>0}$ and $a_{ij}\in\mathbb{Z}_{\geq 0}$,
for each $h=1,\ldo,m$ there is an index $i_h=1,\ldo,k$ such
 that $b_{i_h}=a_{i_hh}=1$ and $a_{ih}=0$ for all $i\neq i_h$. This
 means that $y_h$ appears only in $x_{i_h}$, with coefficient equal to 1.
\end{dimo}

In order to study the possible subdivisions, 
we need to know how non-extremal relations 
decompose in $\NE(X)$. Since in the general the coefficients of such a 
decomposition are not integral, but rational, we need to 
introduce a new type of 
classes in $\NE(X)$, which will allow us to work with integral decompositions.
We refer to \cite{contr} for a more detailed account on this subject.
\begin{defi}
Let $Z$ be an $n$-dimensional, complete, smooth toric variety, and $\gamma\in\NE(Z)$. We say that $\gamma$ is \emph{contractible} if the two
following conditions hold: 
\begin{enumerate}[(i)]
\item
$\gamma$ is primitive, i.~e.\ there exists a primitive collection $P$ 
in $\f_Z$ such that $\gamma=r(P)$. Therefore we have 
\[ \gamma\colon \quad \som{h}=a_1y_1+\cdots+a_ky_k, \] 
with $k\geq 0$, $a_i\in\Z_{> 0}$ for all $i$ and $\langle 
y_1,\ldo,y_k\rangle\in\f_Z$. 
\item if $\nu=\langle z_1,\ldo,z_t\rangle$ is such that $\langle
y_1,\ldo,y_k\rangle +\nu\in\f_Z$ and
$\{z_1,\ld,z_t\}
\cap\{ x_1,\ld,x_h,y_1,\ld,y_k\}=\emptyset$, 
 then
\[ \langle x_1,\ld,\check{x}_i,\ld,x_h,y_1,\ld,y_k\rangle+\nu\in\f_Z
\quad\text{ for all }i=1,\ldo,h.\] 
\end{enumerate} 
\end{defi} 
Comparing with theorem~\ref{reidextr}, we see that this is exactly Reid's description of the geometry of the fan around the walls corresponding to an extremal class.
This description is actually equivalent to the existence of an equivariant 
morphism $\ph_{\gamma}\colon Z\rightarrow Z_{\gamma}$ with connected fibers
 such that for all curves $C\subset Z$
\[ \ph_{\gamma}(C)=\{pt\} \quad \Longleftrightarrow\quad [C]\in\Q_{\,\geq 0}\gamma. \]
In projective varieties, contractible classes have the following interesting 
property:
\begin{teo}[\cite{contr}]
\label{bigteo}
Let $Z$ be an $n$-dimensional, smooth,  projective toric variety. Then
for every $\eta\in\A(Z)\cap\NE(Z)$ there is a decomposition
\[ \eta=m_1 \gamma_1+\cdots+m_r \gamma_r \]
with $\gamma_i$ contractible and $m_i\in\Z_{>0}$ for all $i=1,\ldo,r$.
\end{teo}

\begin{lemma}
\label{gradodue}
Let $X$ be a Fano 4-fold and $\{x_1,x_2,x_3\}$ a primitive collection
in $X$ with relation $x_1+x_2+x_3=x$. Suppose that the relation is not
contractible. Then there are three possible kinds 
of decomposition as a sum of
primitive relations of degree 1:
\begin{enumerate}[(A)]
\item 
$x_1+x=y$, $y+x_2+x_3=2x$. Then the cones $\langle x,x_2,x_3\rangle$,
$\langle x,x_2,y\rangle$, $\langle x,x_3,y\rangle$, $\langle
y,x_1,x_2\rangle$, $\langle y,x_1,x_3\rangle$ are in $\fx$.
\item 
$x_1+w=z$, $z+x_2+x_3=w+x$. Then the cones $\langle
w,x,x_2,x_3\rangle$, $\langle w,x,z,x_2\rangle$, $\langle
w,x,z,x_3\rangle$, $\langle z,x,x_1,x_2\rangle$, $\langle
z,x,x_1,x_3\rangle$ are in $\fx$. 
\item
$x_1+z+w=2x$, $x+x_2+x_3=z+w$. In this case, the cone
$\langle z,w\rangle$ is in $\fx$ and crosses $\langle x,x_2,x_3\rangle$.
\end{enumerate}
\end{lemma}

\vspace{15pt}

\hspace{5pt}(A) \psfig{figure=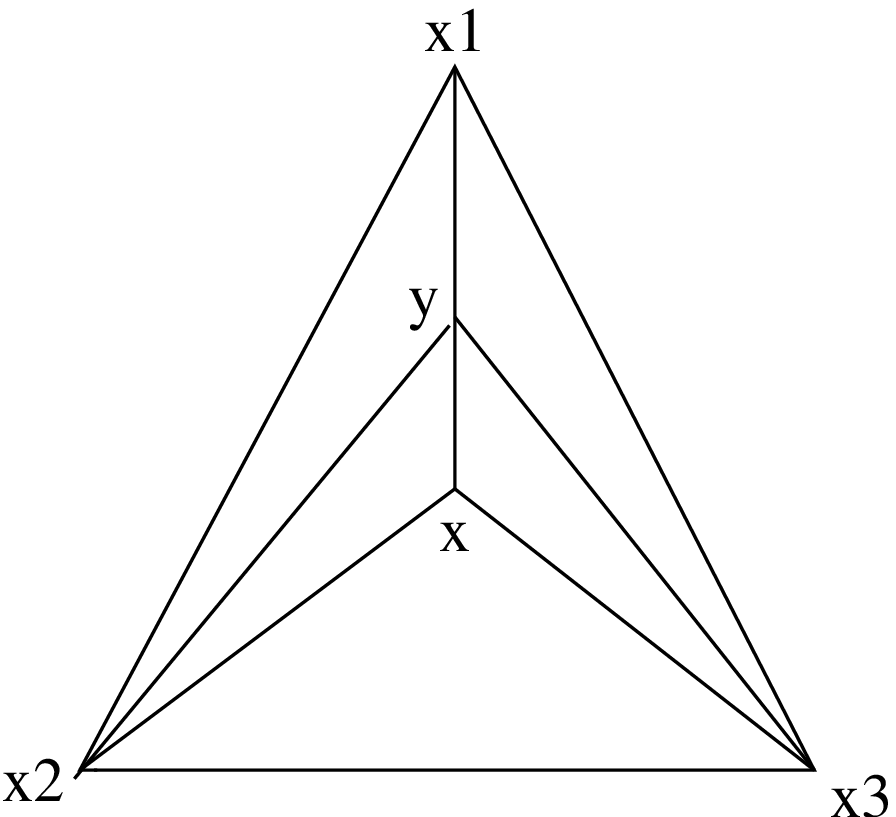,width=4cm}\hspace{20pt}(B)\psfig{figure=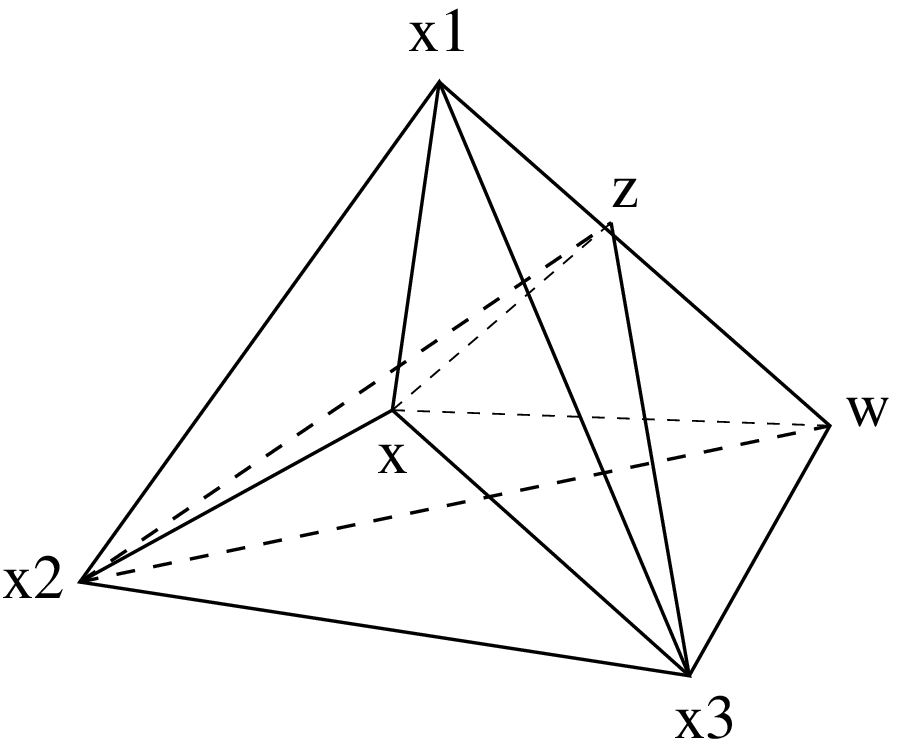,width=5cm}

\vspace{15pt} 

\begin{dimo}
Since the relation $x_1+x_2+x_3=x$ has degree two and it is not contractible,
by theorem~\ref{bigteo} it must be the sum of two relations of degree one.
There are four types of primitive relations of degree one in a toric
Fano 4-fold:
\begin{enumerate}[a.]
\item $y_1+y_2=y_3$;
\item $y_1+y_2+y_3=2y_4$;
\item $y_1+y_2+y_3=y_4+y_5$; 
\item $y_1+y_2+y_3+y_4=3y_5$.
\end{enumerate}
We have to examine all possible combinations.
\begin{description}
\item (a+a): $y_1+y_2-y_3+w_1+w_2-w_3=x_1+x_2+x_3-x$. Then we have
  $y_3=w_1$, $w_3=x$, $y_1=x_1$, $y_2=x_2$ which implies $\langle
  x_1,x_2\rangle\not\in\fx$, a contradiction.
\item (a+b): $y_1+y_2-y_3+w_1+w_2+w_3-2w_4=x_1+x_2+x_3-x$. Then we have
  $y_1=w_4=x$, $y_3=w_1$, $y_2=x_1$, $w_2=x_2$, $w_3=x_3$; this gives
  (A).
\item (a+c):  $y_1+y_2-y_3+w_1+w_2+w_3-w_4-w_5=x_1+x_2+x_3-x$. Then we have
  $w_5=x$, $y_1=w_4$, $y_3=w_1$, $y_2=x_1$, $w_2=x_2$, $w_3=x_3$; this
  gives (B).
\item (a+d): $y_1+y_2-y_3+w_1+w_2+w_3+w_4-3w_5=x_1+x_2+x_3-x$,
  impossible.
\item (b+b):  $y_1+y_2+y_3-2y_4+w_1+w_2+w_3-2w_4=x_1+x_2+x_3-x$,
  impossible.
\item (b+c):   $y_1+y_2+y_3-2y_4+w_1+w_2+w_3-w_4-w_5=x_1+x_2+x_3-x$.
Then we have
  $y_4=w_1=x$, $y_1=w_4$, $y_2=w_5$, $y_3=x_1$, $w_2=x_2$, $w_3=x_3$; 
this gives (C).
\item (b+d): $y_1+y_2+y_3-2y_4+w_1+w_2+w_3+w_4-3w_5=x_1+x_2+x_3-x$,
  impossible.
\item (c+c):
  $y_1+y_2+y_3-y_4-y_5+w_1+w_2+w_3-w_4-w_5=x_1+x_2+x_3-x$. Then we
  have $y_5=x$, $y_4=w_1$, $w_4=y_1$, $w_5=y_2$, $y_3=x_1$, $w_2=x_2$,
  $w_3=x_3$.
We get the relations: $y_1+y_2+x_1=y_4+x$ and $y_4+x_2+x_3=y_1+y_2$. 
Since they both have degree 1, they are extremal, so the cones
$\langle y_1, y_2, y_4,x\rangle$, 
$\langle y_1, y_2, y_4,x_2\rangle$ and $\langle y_1, y_2, y_4,x_3\rangle$
are in $\fx$, which is
impossible.
\item (c+d): $y_1+y_2+y_3-y_4-y_5+w_1+w_2+w_3+w_4-3w_5=x_1+x_2+x_3-x$,
  impossible.
\item (d+d): $y_1+y_2+y_3+y_4-3y_5+w_1+w_2+w_3+w_4-3w_5=x_1+x_2+x_3-x$,
  impossible.
\end{description}
\end{dimo}
\begin{dimo}[Proof of proposition \ref{subdiv}]
Consider the unique cone $\tilde{\sigma}\in\fx$ such that $y_1+y_2+y_3+y_4 \in\RelInt\tilde{\sigma}$, where $\sigma=\langle y_1,y_2,y_3,y_4\rangle$. The generators of $ \tilde{\sigma}$ are described in lemma~\ref{centro}.

(I) If $\dim\tilde{\sigma}=4$, then $\tilde{\sigma}=\sigma$; this is (1).

(II) If $\dim\tilde{\sigma}=3$, 
by lemma \ref{centro} the only
 possibility is $x_1=y_1+y_2\in G(\fx)$
and  $\tilde{\sigma}=\langle x_1,y_3,y_4\rangle$. In $\fx$
$y_1+y_2=x_1$ is a primitive relation of degree 1, so it is
 extremal; $\tilde{\sigma}\in\fx$ implies $\langle
x_1,y_1,y_3,y_4\rangle\in\fx$,  $\langle
x,y_2,y_3,y_4\rangle\in\fx$; this gives (2).

(III) Suppose $\dim\tilde{\sigma}=2$. By lemma \ref{centro}, there are two
possibilities for $\tilde{\sigma}$: $\tilde{\sigma}=\langle
y_1+y_2,y_3+y_4 \rangle$ or $\tilde{\sigma}=\langle
y_1+y_2+y_3,y_4 \rangle$.

(III.1) Let $\tilde{\sigma}=\langle x_1,x_2 \rangle$ with $x_1=y_1+y_2$,
$x_2=y_3+y_4$. Then $y_1+y_2=x_1$ and $y_3+y_4=x_2$ are primitive
relations of degree 1, hence they are extremal. So $\tilde{\sigma}\in\fx$
implies that the cones $\langle x_1,x_2, y_1,y_3 \rangle$, $\langle
x_1,x_2,y_1,y_4 \rangle$, $\langle x_1,x_2,y_2,y_3 \rangle$, $\langle
x_1,x_2,y_2,y_4 \rangle$ are in $\fx$ and this gives (5).

(III.2) Let $\tilde{\sigma}=\langle x_1,y_4 \rangle$ with
$x_1=y_1+y_2+y_3$. Then 
$\langle y_1,y_2,y_3\rangle\not\in\fx$, so $\{y_1,y_2,y_3\}$ contains
a primitive collection. 

(III.2.1) Suppose $\{y_1,y_2,y_3\}$ is primitive: then the primitive
relation is $y_1+y_2+y_3=x_1$. 

(III.2.1.1) If the relation $y_1+y_2+y_3=x_1$ is contractible,
$\tilde{\sigma}\in\fx$ implies that 
$\langle x_1,y_1,y_2,y_4 \rangle$, $\langle x_1,y_1,y_3,y_4 \rangle$ and
$\langle x_1,y_2,y_3,y_4 \rangle$ are in $\fx$ and we get (3).

(III.2.1.2) If $\{y_1,y_2,y_3\}$ is primitive but $y_1+y_2+y_3=x_1$ is not
contractible, then it must be the sum of two  relations of
degree 1, and the possible cases are described in lemma
\ref{gradodue}. We remark that the case (C) never occurs, because
$\fx$ is a subdivision of $\fy$.

(III.2.1.2.1) If the sum is like in (A), then we have two extremal
relations: 
$y_1+x_1=x_2$ and $y_2+y_3+x_2=2x_1$. Since $\tilde{\sigma}=\langle
x_1,y_4\rangle$ is in $\fx$, we get that the cones $\langle x_1, y_2,
y_3,y_4\rangle$, $\langle x_1, x_2, y_2,y_4\rangle$, $\langle x_1,
x_2, y_3,y_4\rangle$, $\langle x_2, y_1, y_2,y_4\rangle$, $\langle
x_2, y_1, y_3,y_4\rangle$ are in $\fx$ too. This is (6).

(III.2.1.2.2) If the sum is like in (B), we have the extremal relations
$y_1+x_3=x_2$ and $x_2+y_2+y_3=x_1+x_3$. Examining the 4-dimensional
cones in $\fx$ given by these relations, we see that the 2-dimensional
cones of $\fx$ contained in $\sigma$ and containing $x_1$ are $\langle
x_1,y_1\rangle$, $\langle x_1,y_2\rangle$, $\langle x_1,y_3\rangle$, $\langle
x_1,x_2\rangle$ and $\langle x_1,x_3\rangle$. Since $\tilde{\sigma}=\langle
x_1,y_4\rangle\in\fx$, it must be $y_4=x_2$ or $y_4=x_3$. If
$y_4=x_2$, then $x_3=y_4-y_1\not\in\sigma$, which is impossible; hence
$y_4=x_3$ and $x_2=y_1+y_4$, and we get (7).

(III.2.2) If $\{y_1,y_2,y_3\}$ is not primitive, it must contain a primitive
collection of order 2: 
we can suppose that it is $\{y_1,y_2\}$. The primitive relation can
not be $y_1+y_2=0$, because $y_1$ and $y_2$ both belong to
$\sigma$. Therefore we get an extremal relation
$y_1+y_2=x_2$. Moreover, since $y_3+x_2=x_1$, $\{y_3,x_2\}$ is a
primitive collection of degree 1, which gives another extremal relation. 
Hence
$\tilde{\sigma}=\langle x_1,y_4 \rangle\in\fx$ implies that
$\langle x_1,x_2,y_1,y_4 \rangle$ and $\langle x_1,x_2,y_2,y_4
\rangle$ are in $\fx$. We claim that also $\langle x_1,y_1,y_3,y_4
\rangle$ and $\langle x_1,y_2,y_3,y_4 \rangle$ are in $\fx$, so the
subdivision is (8). Indeed,
suppose there exists $z\in G(\fx)$, $z\in\langle x_1,y_1,y_3,y_4
\rangle$, different from the vertexs. Then $\{z,x_2\}$ should be a
primitive collection, which is impossible, because it could not have
positive degree. 

(IV) Finally, suppose  $\dim\tilde{\sigma}=1$, i.~e.\
$x_1=y_1+y_2+y_3+y_4\in G(\fx)$. 

(IV.1) We first consider the case where  $\{y_1,y_2,y_3,y_4\}$ is primitive.

(IV.1.1) If $\{y_1,y_2,y_3,y_4\}$ is primitive and 
the relation $y_1+y_2+y_3+y_4=x_1$ is contractible,  then $\langle
x_1,y_1,y_2,y_3\rangle$, $\langle x_1,y_1,y_2,y_4\rangle$,
$\langle x_1,y_1,y_3,y_4\rangle$,
$\langle x_1,y_2,y_3,y_4\rangle$ are in $\fx$, so we get (4). 

(IV.1.2) If $\{y_1,y_2,y_3,y_4\}$ is primitive and 
$y_1+y_2+y_3+y_4=x_1$ is not contractible, then by
lemma \ref{puh}
it must be $Y\simeq\pr{4}$: indeed, all the 3-dimensional faces of
$\sigma$ are in $\fx$, so there must be at least one generator 
of $\fx$ different from $x_1$ in the interior
of $\sigma$.
Again by lemma~\ref{puh}, we know that
$\fx$ has exactly two
generators in the interior of $\sigma$, $x_1$ and $x_2$, with 
 primitive relations $(-x_1)+x_2=y_1$ and
$x_1+y_1=x_2$.
We get two effective relations, $x_1+y_1=x_2$ and $y_2+y_3+y_4+x_2=2x_1$, whose
sum is  $y_1+y_2+y_3+y_4=x_1$. Since there are no other new generators
in $\sigma$, these relations must be  both contractible, therefore the
cones $\langle x_1,x_2,y_2,y_3\rangle$,
$\langle x_1,x_2,y_2,y_4\rangle$,
$\langle x_1,x_2,y_3,y_4\rangle$,
$\langle x_1,y_2,y_3,y_4\rangle$,
$\langle x_2,y_1,y_2,y_3\rangle$,
$\langle x_2,y_1,y_2,y_4\rangle$,
$\langle x_2,y_1,y_3,y_4\rangle$ are in $\fx$ and we get (14).

(IV.2) Suppose now that $\{y_1,y_2,y_3\}$ is a primitive collection. The
primitive relation can not be $y_1+y_2+y_3=z_1+z_2$, because the cone
$\langle z_1,z_2\rangle$ would cross $\langle y_1,y_2,y_3\rangle$ and
$\fx$ would not be a subdivision of $\fy$. Neither can the relation be 
$y_1+y_2+y_3=2x_2$, because $\{y_1,y_2,y_3\}$ is a part of a basis of
the lattice. Thefore the primitive relation is $y_1+y_2+y_3=x_2$. Then
$\{x_2,y_4\}$ is a primitive collection with relation $x_2+y_4=x_1$,
which is extremal.

(IV.2.1) If $y_1+y_2+y_3=x_2$ is contractible, then the subdivision is given by
the cones $\langle x_1,x_2,y_1,y_2\rangle$,
$\langle x_1,x_2,y_1,y_3\rangle$,
$\langle x_1,x_2,y_2,y_3\rangle$,
$\langle x_1,y_1,y_2,y_4\rangle$,
$\langle x_1,y_1,y_3,y_4\rangle$,
$\langle x_1,y_2,y_3,y_4\rangle$, and this gives (9).

(IV.2.2) If $y_1+y_2+y_3=x_2$ is not contractible, then it is sum of two
primitive relations of degree 1, and we have to examine the cases
given by lemma~\ref{gradodue};
 we remark again that case (C) is not possible here,
because $\fx$ is a subdivision of $\fy$.

(IV.2.2.1) If the sum is like in (A), then we have primitive relations
$y_1+x_2=x_3$ and $x_3+y_2+y_3=2x_2$. We claim that this is
impossible. Indeed, by the extremality of $x_2+y_4=x_1$ and
$y_1+x_2=x_3$, we get that $\langle x_1,y_1\rangle\in\fx$; on the
other hand, we have $y_4+x_3=x_1+y_1$, so $\{y_4,x_3\}$ would be a primitive
collection of degree zero. 

(IV.2.2.2) We have two possible  decompositions  as in case (B) of lemma
\ref{gradodue}: 
\begin{description}
\item
(B1) $y_1+x_3=x_1$, $y_2+y_3+x_1=x_2+x_3$. Then the cones  $\langle
x_1,y_1,y_2,y_4\rangle$,  $\langle x_1,y_1,y_3,y_4\rangle$,  
$\langle x_1,x_2,y_1,y_2\rangle$,  
$\langle x_1,x_2,y_1,y_3\rangle$,  $\langle
x_1,x_3,y_2,y_4\rangle$,
\linebreak  $\langle x_1,x_3,y_3,y_4\rangle$,  $\langle
x_2,x_3,y_2,y_3\rangle$,  $\langle x_1,x_2,x_3,y_2\rangle$,  $\langle
x_1,x_2,x_3,y_3\rangle$ are in
$\fx$, and we get (10).  

\item
(B2) $y_1+x_1=x_3$, $y_2+y_3+x_3=x_1+x_2$ . Then the cones
$\langle x_1,y_2,y_3,y_4\rangle$, $\langle x_3,y_1,y_2,y_4\rangle$,  $\langle
x_3,y_1,y_3,y_4\rangle$,  $\langle x_1,x_2,y_2,y_3\rangle$, $\langle
x_1,x_3,y_2,y_4\rangle$, 
\linebreak $\langle
x_1,x_3,y_3,y_4\rangle$, $\langle x_2,x_3,y_1,y_2\rangle$,  $\langle
x_2,x_3,y_1,y_3\rangle$, $\langle x_1,x_2,x_3,y_2\rangle$, 
\linebreak $\langle x_1,x_2,x_3,y_3\rangle$
are in $\fx$, and we get (15).
\end{description}
The last subdivision can occur only if $Y\simeq\pr{4}$,
because there are two generators of $\fx$ in the interior of
$\sigma$. 

(IV.3) Finally, suppose $\{y_1,y_2,y_3,y_4\}$ does not contain primitive
relations of order 3. Then we can assume $\{y_1,y_2\}$ primitive, with
relation $y_1+y_2=x_2$. Since $y_3+y_4+x_2=x_1$, the cone $\langle
y_3,y_4,x_2\rangle$ can not be in $\fx$, so $\{y_3,y_4,x_2\}$ contains a
primitive collection. 

(IV.3.1) If $\{y_3,y_4\}$ is primitive, with relation $y_3+y_4=x_3$, then
$x_2+x_3=x_1$ is also a primitive relation.    
Hence the subdivision is given by the cones $\langle
x_1,x_2,y_1,y_3\rangle$,   $\langle x_1,x_2,y_2,y_3\rangle$,  $\langle
x_1,x_2,y_1,y_4\rangle$,  $\langle x_1,x_2,y_2,y_4\rangle$,  $\langle
x_1,x_3,y_1,y_3\rangle$,  $\langle x_1,x_3,y_1,y_4\rangle$,  $\langle
x_1,x_3,y_2,y_3\rangle$,  $\langle x_1,x_3,y_2,y_4\rangle$, and we get
(12).
 
(IV.3.2) If  $\{y_3,x_2\}$ is primitive, with relation $y_3+x_2=x_3$, then
$y_4+x_3=x_1$ is also a primitive relation.   Hence the subdivision is
given by the cones
 $\langle x_1,y_1,y_3,y_4\rangle$,  $\langle x_1,y_2,y_3,y_4\rangle$,
 $\langle x_1,x_2,y_1,y_4\rangle$,  $\langle x_1,x_2,y_2,y_4\rangle$,
 $\langle x_1,x_3,y_1,y_3\rangle$,  $\langle x_1,x_3,y_2,y_3\rangle$,
 $\langle x_1,x_2,x_3,y_1\rangle$,  $\langle x_1,x_2,x_3,y_2\rangle$;
 this gives (13).

(IV.3.3) The last possibility is that  $\{y_3,y_4,x_2\}$ is a primitive
collection. 

(IV.3.3.1) If $y_3+y_4+x_2=x_1$ is contractible, then
the cones  $\langle x_1,y_1,y_3,y_4\rangle$,  $\langle
x_1,y_2,y_3,y_4\rangle$,  $\langle x_1,x_2,y_1,y_3\rangle$,  $\langle
x_1,x_2,y_2,y_3\rangle$,  $\langle x_1,x_2,y_1,y_4\rangle$,  $\langle
x_1,x_2,y_2,y_4\rangle$ are in $\fx$ and this gives (11).

(IV.3.3.2) If $\{y_3,y_4,x_2\}$ is primitive but $y_3+y_4+x_2=x_1$ is not
contractible,  then it must decompose as a sum of two primitive
relations of degree 1 as described in lemma
\ref{gradodue}.

We have two possible  decompositions  as in case (A) of lemma
\ref{gradodue}: 
\begin{description}
\item
(A1) $y_3+x_1=x_3$, $y_4+x_2+x_3=2x_1$. Then the cones  $\langle
x_1,x_2,y_1,y_4\rangle$,  $\langle x_1,x_2,y_2,y_4\rangle$,  $\langle
x_3,y_1,y_3,y_4\rangle$,  $\langle x_3,y_2,y_3,y_4\rangle$,  $\langle
x_1,x_3,y_1,y_4\rangle$,  
\linebreak $\langle x_1,x_3,y_2,y_4\rangle$,  $\langle
x_2,x_3,y_1,y_3\rangle$,  $\langle x_2,x_3,y_2,y_3\rangle$,  $\langle
x_1,x_2,x_3,y_1\rangle$,  
\linebreak $\langle x_1,x_2,x_3,y_2\rangle$ are in
$\fx$, and we get (16).  

\item
(A2) $x_1+x_2=x_3$, $y_3+y_4+x_3=2x_1$. Then the cones  $\langle
x_1,y_1,y_3,y_4\rangle$,  $\langle x_1,y_2,y_3,y_4\rangle$,  $\langle
x_1,x_3,y_1,y_3\rangle$,  $\langle x_2,x_3,y_1,y_3\rangle$,  $\langle
x_1,x_3,y_2,y_3\rangle$,  
\linebreak $\langle x_2,x_3,y_2,y_3\rangle$,  $\langle
x_1,x_3,y_1,y_4\rangle$,  $\langle x_2,x_3,y_1,y_4\rangle$,  $\langle
x_1,x_3,y_2,y_4\rangle$,  
\linebreak $\langle x_2,x_3,y_2,y_4\rangle$ are in
$\fx$, and we get (17).
\end{description}
We remark that both subdivisions can occur only when $Y\simeq\pr{4}$,
because there are two generators of $\fx$ in the interior of
$\sigma$. 

We claim that $y_3+y_4+x_2=x_1$ can not decompose as in case (B) of
lemma \ref{gradodue}. Indeed, either the decomposition is given by
$y_3+w=z$ and $y_4+x_2+z=w+x_1$, or by $x_2+w=z$ and
$y_3+y_4+z=w+x_1$; in both cases we get $\langle
z,y_3,x_1,x_2\rangle\in\fx$. Now, $y_1+y_2=x_2$ is extremal, so it
must be $z=y_1$ or $z=y_2$. This implies $w\not\in\sigma$, a contradiction.
\end{dimo}

\begin{dimo}[Proof of theorem \ref{facto}]
We know by proposition \ref{subdiv} that the subdivision of every
4-dimensional cone of $\fy$ is given by a sequence of
star-subdivisions.  We just have to check that 
we can order these star-subdivisions in a way which is compatible with
the orders in all 4-dimensional cones of $\fy$.

We order the cones corresponding to centers of star-subdivisions in
such a way that: 
\begin{enumerate}[(i)]
\item the dimensions of the cones are non-increasing 
\item each cone appears in the fan obtained with the preceeding 
star-subdivi- \newline sions.
\end{enumerate}  
It is easy to see directly, looking at the subdivisions given by
proposition~\ref{subdiv}, that this order works. 
More precisely, in all the sequences of
star-subdivisions given by
proposition~\ref{subdiv}, except for (5) and (10),
there is a unique cone of maximal dimension. In (5) and (10) both
cones of maximal dimension can be chosen for the first star-subdivision.
Then we can use induction.
\end{dimo}

\medskip

\small
\label{ultima}

\bigskip

\noindent Cinzia Casagrande \\
Dipartimento di Matematica\\
Universit\`a di Roma ``La Sapienza''\\
piazzale Aldo Moro, 5\\
00185 Roma

\noindent ccasagra@mat.uniroma1.it
\end{document}